\documentclass[a4paper,11pt]{article}

\usepackage{bm} \usepackage{amsmath} \usepackage{amsfonts}
\usepackage{xcolor} \usepackage{graphicx} \usepackage{inputenc}
\usepackage[top = 0.5in, left = 0.5 in, right = 0.5 in, bottom = 0.8 in]{geometry}
\usepackage[section]{placeins}
\usepackage{float}
\usepackage{array}
\usepackage{datetime}
\usepackage{enumerate}
\usepackage{nomencl} 
\makenomenclature

\usepackage{amssymb}
\usepackage{mathbbol}
\usepackage{bbm}
\usepackage{accents}
\usepackage{amsthm}
\usepackage{color}
\usepackage{graphics}
\usepackage{epsfig}
\usepackage{cancel}
\usepackage{epstopdf}
\usepackage{hyperref}
\usepackage{xcolor}
\usepackage{upgreek}
\usepackage{caption}
\usepackage{multirow}
\usepackage{subcaption}
\usepackage{setspace}
\usepackage{comment}

\hypersetup{
colorlinks,
linkcolor={blue},
citecolor={blue},
urlcolor={blue},
}
\usepackage[numbers,sort&compress]{natbib}
\usepackage{stmaryrd}

\usepackage{mdframed} 
\newmdtheoremenv[linecolor=white,leftmargin=0,rightmargin=0,backgroundcolor=white!10,innertopmargin=0pt,ntheorem]{sidenote}{Remark}[section]

\usepackage{array}
\newcolumntype{L}[1]{>{\raggedright\let\newline\\\arraybackslash\hspace{0pt}}m{#1}}
\newcolumntype{C}[1]{>{\centering\let\newline\\\arraybackslash\hspace{0pt}}m{#1}}
\newcolumntype{R}[1]{>{\raggedleft\let\newline\\\arraybackslash\hspace{0pt}}m{#1}}

\usepackage{tabu}

\numberwithin{equation}{section}

\title{Plane stress finite element modelling of arbitrary compressible hyperelastic materials}
\author{Masoud Ahmadi\textsuperscript{1}, Andrew McBride\textsuperscript{1}, Paul Steinmann\textsuperscript{1,2}, Prashant Saxena\textsuperscript{1}\thanks{Corresponding author email: prashant.saxena@glasgow.ac.uk} \\[2ex]
\textsuperscript{1}{\small Glasgow Computational Engineering Centre, James Watt School of Engineering}\\
{\small University of Glasgow, Glasgow G12 8LT, UK} \\
\textsuperscript{2}{\small Institute of Applied Mechanics, Friedrich-Alexander University Erlangen-Nürnberg, D-91052, Erlangen, Germany}}

\date{}

\allowdisplaybreaks 

\begin{document}


\maketitle


\begin{abstract}
Modelling the large deformation of hyperelastic solids under plane stress conditions for arbitrary compressible and nearly incompressible material models is challenging.
This is in contrast to the case of full incompressibility where the out-of-plane deformation can be entirely characterised by the in-plane components.
A rigorous general procedure for the incorporation of the plane stress condition for the compressible case (including the nearly incompressible case) is provided here, accompanied by a robust and open source finite element code.
An isochoric/volumetric decomposition is adopted for nearly incompressible materials yielding a robust single-field finite element formulation.
The nonlinear equation for the out-of-plane component of the deformation gradient is solved using a Newton--Raphson procedure nested at the quadrature point level. 
The model's performance and accuracy are made clear via a series of simulations of benchmark problems. 
Additional challenging numerical examples of composites reinforced with particles and fibres further demonstrate the capability of this general computational framework. 

\end{abstract}


\section{Introduction}\label{sec.int}

Soft materials, prevalent in elastomers, biological tissues, and engineering applications, exhibit nonlinear hyperelastic behaviour under large deformations \citep{yoda1998elastomers, shahzad2015mechanical, gultekin2019quasi,  hossain2013more}.
Planar approximations (plane stress and plane strain) are appropriate and computationally beneficial for specific scenarios; particularly, the plane stress approximation becomes crucial when the structure is thin.
However, details of the general finite element framework to impose the plane stress condition for arbitrary hyperelastic material models at large deformations are sparse for the compressible (and nearly incompressible) regime.

In the past decades, significant advancements have been made in the field of nonlinear elasticity enabling the study of the complex nature of hyperelastic materials and their composites.
Computational methods, specifically the finite element method, have been widely employed to study the mechanical response of these materials under large deformations \citep{reese2002equivalent, liu2024computational, doghri2013mechanics, bijelonja2005finite, li2000numerical, steinmann1994micropolar},
providing valuable insights for engineers and researchers.
Research on finite element modelling of nonlinear elasticity tends to focus on the three-dimensional problem or the plane strain approximation thereof
\citep{shojaei2019compatible, angoshtari2017compatible, auricchio2013approximation, brink1996some, chavan2007locking}. 
In addition, truly two-dimensional elasticity formulations are employed, which, inspired by Edwin Abbott's 1884 novel \textit{Flatland: A Romance in Many Dimensions}, will henceforth be termed ``flatland'', a term also adopted by others \citep[see, e.g.,][]{Rizzo2003}.
In the flatland approach, all kinetic and kinematics quantities (stress and strain) are completely restricted to two dimensions. 
Despite the large body of work on three-dimensional, plane strain and flatland problems, there are very few algorithmic expositions on the imposition of the plane stress assumption when modelling general classes of compressible hyperelastic materials under large deformation.
This deficiency is addressed here. 


The general procedure for incorporating plane stress conditions has been discussed in textbooks \citep{zienkiewicz2005finite}. 
The plane stress assumption is routinely made for structural elements such as plates and shells. 
For example, \citet{viebahn2017simple} developed a shell theory for compressible thin hyperelastic shells wherein the out-of-plane component of the deformation was determined approximately from the solution of the nonlinear plane stress constraint. 
\citet{klinkel2002using} developed algorithms to enforce zero-stress conditions in structural elements such as beams and shells. Their approach condenses arbitrary material laws to account for these stress conditions at each integration point on the element level.
The distinct kinematics of shell theory hinders the direct translation of this approach to continuum elements. 
For further algorithmic details on the imposition of the plane stress condition for fully nonlinear structural elements, see \citet{korelc2016automation} where the finite element framework AceGEN \citep{korelc2009automation} is used for implementation. 
\citet{steinmann1997fe} proposed a formulation for hyperelastic materials wherein the plane stress condition is satisfied weakly, that is, it is not enforced directly at the level of the quadrature point.
The basis of the work presented in this manuscript is to develop a procedure for the strong imposition of the plane stress assumption for general classes of compressible hyperelastic materials that build directly upon and generalise the contribution of \citet{pascon2019large}. 
Their work exploited the structure of specific classes of material models when constructing the plane stress relation. 
However, their approach is not directly applicable to general constitutive models, including those important ones that impose an isochoric/volumetric decomposition \citep{miehe1994aspects}. 
The present work extends these works, providing a robust finite element algorithm and accompanying code for modelling complex material behaviour.

Given the limited scope of prior investigations into plane stress formulations for compressible (and nearly incompressible) hyperelastic materials, this study presents a general framework for their nonlinear finite element analysis, accompanied by an open source code \citep{madeal2024}.
This approach accommodates arbitrary hyperelastic strain energy functions, including those with an isochoric/volumetric split.
The framework is validated using challenging nonlinear problems, including complex cases involving reinforced stiff particles and fibres.
Furthermore, comparative analyses evaluate different approaches considering aspects of three dimensions, and planar approximations thereof as well as flatland modelling.
In passing we note that a three-field mixed formulation based on the Hu--Washizu variational principle \citep{hu1984variational, simo1986variational, wriggers2008nonlinear} is adopted to address the issue of volumetric locking in hyperelastic materials for non-plane-stress cases. 

The remainder of this contribution is structured as follows.
A brief overview of nonlinear hyperelasticity is provided in Section \ref{sec: math}.
Section \ref{sec: 2d} presents two-dimensional modelling, including the planar approximation of plane stress (and for completeness plane strain) as well as the flatland approach.
To validate the finite element models, two benchmark problems are examined in Section \ref{sec: numerical}, comparing the performance of three- and two-dimensional approaches. 
Additional challenging numerical examples of inhomogeneous problems, i.e., composites reinforced by particles and fibres, are explored to assess the plane stress (and plane strain) condition for both compressible and nearly incompressible materials.
Finally, Section \ref{sec: conclusion} draws some conclusions.

\paragraph*{Mathematical notation:}

Both direct and index notations are used for algebra over vectors and tensors.
In the direct notation, the scalars are written in italic font, vectors and 2nd-order tensors are written in bold font, and 4th-order tensors are denoted by blackboard letters such as $\mathbb{D}$.
Einstein’s summation convention is used to sum over repeated indices, that is, $a_i b_i = \mathbf{a} \cdot \mathbf{b} = \sum\limits_{i=1}^n a_i b_i$, where $n$ is the dimension of $a$ and $b$.
The inner product between two 2nd-order tensors is defined as $\mathbf{A}:\mathbf{B} = A_{ij}B_{ij}$.
When a second or higher-order tensor is involved in a single contraction product, the dot symbol ($\cdot$) is omitted, following the convention commonly used in the literature.
The 2nd-order identity tensor $\mathbf{I}$ has components $I_{ij} = \delta_{ij}$, where $\delta_{ij}$ is the Kronecker delta.
The fourth-order symmetric identity tensor $\mathbb{I}$ is defined as $I_{ijkl}=\frac{1}{2}\left[\delta_{ik}\delta_{jl}+\delta_{il}\delta_{jk}\right]$ such that $\mathbb{I} \colon \mathbf{A}= \text{sym}(\mathbf{A})$.
The dyadic-type tensor products $\otimes$ and $\odot$ for tensors $\mathbf{A}$ and $\mathbf{B}$ are defined by
\begin{align}
    [{\mathbf{A}} \otimes {\mathbf{B}}]_{ijkl} = A_{ij} B_{kl} 
    && \text{and} &&
    [{\mathbf{A}} \odot {\mathbf{B}}]_{ijkl} = \frac{1}{2} \left[ A_{ik} B_{jl} + A_{il} B_{jk} \right] \, .
\end{align}

The Lagrangian and Eulerian gradient of a tensor ($\bullet$) are respectively represented by $\bm{\nabla}_X \, (\bullet) = \dfrac{\partial\left( \bullet \right)}{\partial X_I}\otimes \mathbf{E}_I$ and $\bm{\nabla}_x \, (\bullet) = \dfrac{\partial\left( \bullet \right)}{\partial x_i}\otimes\mathbf{e}_i \,$. 
Corresponding notation is applied to the divergence operator.

\section{A brief on hyperelasticity}\label{sec: math}

This section briefly introduces the essentials for modelling hyperelastic solids under large deformation. 
For detailed expositions on the topic, see standard references such as \citet{ogden1997non} and \citet{holzapfel2000nonlinear}.

\subsection{Kinematics and kinetics}

Consider a hyperelastic solid body occupying a domain $\Omega_X$ in its material configuration with boundary $\Gamma$.
The Lagrangian description of the position of a point in the body is denoted by $\mathbf{X} = X_I \, \mathbf{E}_I \in \Omega_X$ with material basis vectors $\mathbf{E}_I \,$.
Upon a quasi-static deformation, the body occupies the spatial configuration $\Omega_x$ with points with position $\mathbf{x} = x_i \, \mathbf{e}_i$, with spatial basis vectors $\mathbf{e}_i \,$.

The points $\mathbf{X}$ and $\mathbf{x}$ are related by an invertible mapping $\bm{\chi}: \Omega_X \to \Omega_x$ such that $\mathbf{x} = \bm{\chi}(\mathbf{X})$ and the displacement vector follows as $\mathbf{u} (\mathbf{X}) = \mathbf{x} - \mathbf{X}$.
The corresponding deformation gradient, $\mathbf{F}$, is defined by
\begin{equation}
    \mathbf{F} = \frac{\partial \mathbf{x}}{\partial \mathbf{X}} = \bm{\nabla}_X \bm{\chi}.
\end{equation}
The determinant of the deformation gradient is defined by $J = \det{\mathbf{F}}$, where for incompressible materials $J = 1$.
The right Cauchy–Green deformation tensor is defined by $\mathbf{C} = \mathbf{F}^T \, \mathbf{F}$, and the Green--Lagrange
strain tensor
by $\mathbf{E} = \frac{1}{2} \left[ \mathbf{C} - \mathbf{I} \right]$.
The left Cauchy--Green deformation tensor is defined by $\mathbf{b} = \mathbf{F} \, \mathbf{F}^T \,$.

\paragraph{Stress measures:}

The Cauchy stress is denoted as $\bm{\sigma}$.
The Piola stress, $\mathbf{P}$, and the Piola--Kirchhoff stress, $\mathbf{S}$, are related to $\bm{\sigma}$ as
\begin{equation}
    \mathbf{S} = \mathbf{F}^{-1} \mathbf{P} = J \mathbf{F}^{-1} \bm{\sigma} \mathbf{F}^{-T} = \mathbf{S}^T . 
\end{equation}
Finally, the Kirchhoff stress 
is defined as
\begin{equation}
    \bm{\tau} = J \bm{\sigma} = \mathbf{F} \mathbf{S} \mathbf{F}^T .
\end{equation}

\paragraph{Balance of linear momentum:}


The balance of linear momentum for a continuum body (in the absence of body forces) in Lagrangian and Eulerian description are respectively given by 
\begin{align}
    \bm{\nabla}_X\cdot\mathbf{P}=\mathbf{0}
    && \text{and} &&
    \bm{\nabla}_x \cdot \bm{\sigma} = \bm{0} \, .
\end{align}
The details for approximately solving the boundary value problem using a one-field finite element formulation are presented in Appendix \ref{sec: d and l} and for a three-field formulation in Appendix \ref{app: 3field}.


\paragraph{Constitutive relations:}
Isotropic hyperelasticity is characterised by a strain energy density function $\psi$ such that 
\begin{equation}
    \mathbf{S} = \frac{\partial \psi(\mathbf{E})}{\partial \mathbf{E}} \, .
\end{equation}
The deformation can be decomposed into volumetric and isochoric parts, aiding modelling of nearly incompressible materials \citep{holzapfel2000nonlinear}.
The multiplicative decomposition of the deformation gradient into volumetric and isochoric parts is given as
\begin{equation}
    \mathbf{F}=J^{1/3} \, \widehat{\mathbf{F}} \,
    , 
\end{equation}
where $\widehat{\mathbf{F}}$ is isochoric.
Hence, the isochoric part of the right and left Cauchy--Green deformation tensors follow respectively as 
\begin{align}
    \widehat{\mathbf{C}} ={\widehat{\mathbf{F}}^T} \, \widehat{\mathbf{F}}=J^{-2/3} \, \mathbf{F}^T \, \mathbf{F}=J^{-2/3} \, \mathbf{C} && \text{and} && 
    \widehat{\mathbf{b}} =\widehat{\mathbf{F}} \, {\widehat{\mathbf{F}}^T}=J^{-2/3} \, \mathbf{F} \, \mathbf{F}^T=J^{-2/3} \, \mathbf{b} \, .
\end{align}
The strain energy density function may be additively decomposed  based on $J$ and $\widehat{\mathbf{C}}$ as
\begin{align}
\label{eq: psi split}
    \psi=\psi (\widehat{\mathbf{C}},J)=\psi_{\text{iso}}(\widehat{\mathbf{C}})+\psi_{\text{vol}}\left(J\right) ,
\intertext{or, for the case of isotropy, in terms of $J$ and $\widehat{\mathbf{b}}$ as}
    \psi=\psi (\widehat{\mathbf{b}},J) = \psi_{\text{iso}}(\widehat{\mathbf{b}})+\psi_{\text{vol}}\left(J\right) .
\end{align}
The bulk modulus, $\kappa = \lambda + 2/3 \, \mu$, is the sole material parameter that appears in the volumetric part of the energy density function, i.e., $\psi_{\text{vol}} (J) = \kappa \, \mathcal{G} (J)$.
Various formulations for $\mathcal{G} \left(J\right)$ have been proposed, including
\begin{subequations}\label{eq:gfunctions}
\begin{align}
    \mathcal{G} \left(J\right) &= \frac{1}{2} \left[J-1\right]^2, \\ 
    \mathcal{G} \left(J\right) &= \frac{1}{4} \left[J^2-1-2\ln J\right], \label{eq: used G model}\\
    \mathcal{G} \left(J\right) &= \frac{1}{2} \left[\ln J\right]^2, \\
    \mathcal{G} \left(J\right) &= J \ln J-J+1 .
\end{align}
\end{subequations} 
For further discussion on these choices of volumetric strain energy functions, refer to \citep{doll2000development, hartmann2003polyconvexity}.
The decomposition of the strain energy function allows for the corresponding decomposition of the various stress tensors, which are advantageous for numerical implementation, as elaborated upon in Appendix~\ref{app: stress split}.

For isotropic materials, 
the strain energy function can be written as a function of the invariants of $\widehat{\mathbf{C}}$ and $\widehat{\mathbf{b}}$ as
\begin{equation}
    \psi = \psi_{\text{iso}} (I_{\widehat{C}} \, , II_{\widehat{C}}) + \psi_{\text{vol}} (J) ,
\end{equation}
where the scalar invariants are defined by
\begin{subequations}
\begin{align}
    I_{\widehat{C}} &= I_{\widehat{b}} = \text{tr} \, \widehat{\mathbf{C}} = \text{tr} \, \widehat{\mathbf{b}} \, , \\
    II_{\widehat{C}} &= II_{\widehat{b}} = \frac{1}{2} \left[ [\text{tr} \, \widehat{\mathbf{C}}]^2+\text{tr} \, \widehat{\mathbf{C}}^2 \right] = \frac{1}{2} \left[ [\text{tr} \, \widehat{\mathbf{b}}]^2+\text{tr} \, \widehat{\mathbf{b}}^2 \right] \, . 
\end{align}
\end{subequations}
It can be shown that
\begin{equation}\label{eq:hatinavriants}
    I_{\widehat{C}} = J^{-2/3} \, I_C, \quad II_{\widehat{C}} = J^{-4/3} \, {II}_C, 
\end{equation}
where $I_C$ and $II_C$ 
are the scalar invariants of the right Cauchy--Green deformation tensor $\mathbf{C}$.




\section{Two-dimensional approximations and flatland approach}\label{sec: 2d}

In finite element analysis for nonlinear hyperelasticity, the two-dimensional approximation of a three-dimensional problem is often reasonable and brings significant benefits, such as simplified geometry and decreased computational cost.
However, the modelling error associated with this dimension reduction must be appropriately accounted for. 
The following pertinent approaches to two-dimensional modelling are now discussed: plane stress and flatland and, for the sake of completeness, plane strain.

\subsection{Plane stress} 
The plane stress approximation is applied to thin structures by assuming that the stress vanishes along the thickness direction.
The components of the deformation gradient $\mathbf{F}$ and right Cauchy--Green deformation tensor $\mathbf{C}$ become
\begin{align}
    [\mathbf{F}]=\left[\begin{matrix} F_{11}&F_{12}&0 \\ 
            F_{21}&F_{22}&0 \\
            0&0&F_{33} \\
\end{matrix}\right] 
&& \text{and} &&
[\mathbf{C}]= \left[\begin{matrix} F_{11}^2 + F_{21}^2 & F_{11} F_{12} + F_{21} F_{22} & 0 \\ 
    F_{11} F_{12} + F_{21} F_{22} & F_{12}^2 + F_{22}^2 & 0 \\
    0 & 0 & F_{33}^2 \\
    \end{matrix}\right] .  
\end{align}
Accordingly, 
$J = 
F_{33} \, \left[ F_{11} \, F_{22} - F_{12} \, F_{21} \right]$.
Here, determining $C_{33}=F^2_{33}$ is not trivial, 
requiring the use of the plane stress condition which sets the out-of-plane stress component $S_{33}$ to zero.
To avoid ambiguity, the right Cauchy--Green tensor $\mathbf{C} \in \mathbb{R}^{3} \otimes \mathbb{R}^{3}$ is decomposed into in-plane and out-of-plane components as 
\begin{equation}\label{eq: c dec}
    \mathbf{C} = \mathbf{C}_\parallel + \mathbf{C}_\perp,
\end{equation}
where 
\begin{equation}
    [\mathbf{C}_\parallel] = \left[\begin{matrix}     \overline{\mathbf{C}}& \mathbf{0} \\ 
        \mathbf{0}&0 \\
    \end{matrix}\right], \quad
    [\mathbf{C}_\perp] = \left[\begin{matrix}       \overline{\mathbf{0}} &\mathbf{0} \\ 
        \mathbf{0}&C_{33} \\
    \end{matrix}\right], \quad
    \overline{\mathbf{C}}, \overline{\mathbf{0}} \in \mathbb{R}^{2} \otimes \mathbb{R}^{2}, \quad \ C_{33} \in \mathbb{R} \, .
\end{equation}
This decomposition separates $\mathbf{C}$ into the reduced tensor $\overline{\mathbf{C}}$ and the scalar $C_{33}$. 
A similar decomposition of the Piola--Kirchhoff stress $\mathbf{S} \in \mathbb{R}^{3} \otimes \mathbb{R}^{3}$ gives
\begin{equation}
    \mathbf{S} = \mathbf{S}_\parallel + \mathbf{S}_\perp,
\end{equation}
where 
\begin{equation}
    [\mathbf{S}_\parallel] = \left[\begin{matrix}     \overline{\mathbf{S}}&\mathbf{0} \\ 
        \mathbf{0}&0 \\
    \end{matrix}\right], \quad
    [\mathbf{S}_\perp] = \left[\begin{matrix}       \overline{\mathbf{0}}&\mathbf{0} \\ 
        \mathbf{0}&S_{33} \\
    \end{matrix}\right], \quad
    \overline{\mathbf{S}} \in \mathbb{R}^{2} \otimes \mathbb{R}^{2}, \quad S_{33}  \in \mathbb{R} \, .
\end{equation}
The value of $C_{33}$ as a function of $\overline{\mathbf{C}}$ needs to be calculated at every quadrature point within each finite element.
To do so, the nonlinear equation $S_{33} = 0$ is solved using a Newton--Raphson scheme such that 
\begin{equation}
    S_{33}\big (\overline{\mathbf{C}}, C_{33} \big) \big|_{\overline{\mathbf{C}}} + \frac{\partial S_{33} \big (\overline{\mathbf{C}}, C_{33} \big)}{\partial C_{33}}\bigg|_{ \overline{\mathbf{C}}} \ \text{d} C_{33} = 0 ,
\end{equation}
in order to update the solution as
\begin{equation}
    C_{33} \leftarrow C_{33} + \text{d} C_{33},
\end{equation}
until convergence is achieved, resulting in the solution $C_{33} = \acute{C}_{33} (\overline{\mathbf{C}})$.
 
A crucial step is to update the incremental constitutive tensor to account for the plane stress approximation.
The in-plane tensor $\overline{\mathbb{C}} \in  \mathbb{R}^{2} \otimes \mathbb{R}^{2} \otimes \mathbb{R}^{2} \otimes \mathbb{R}^{2}$, is given by
\begin{equation}\label{eq: C_IJKL for pstress}
    \overline{\mathbb{C}} = 2 \, \frac{\partial \overline{\mathbf{S}}}{\partial \overline{\mathbf{C}}} = 2 \left[ \, \frac{\partial \overline{\mathbf{S}}}{\partial \overline{\mathbf{C}}} \bigg|_{\acute{C}_{33}} + \frac{\partial \overline{\mathbf{S}}}{\partial \acute{C}_{33}} \bigg|_{\overline{\mathbf{C}}} \otimes \frac{\partial \acute{C}_{33}}{\partial \overline{\mathbf{C}}} \right] .
\end{equation}
The procedure to compute $\overline{\mathbb{C}}$ for a neo-Hookean strain energy density function is now given.
Note however that the approach is general and can be applied to any hyperelastic model.


\subsubsection{Specialisation to the neo-Hookean model}
\label{sec: neo-Hookean derivations}

The neo-Hookean strain energy density function is given by
\begin{equation}\label{eq:neo-hook function}
    \psi = \psi_{\text{iso}}(I_{\widehat{C}}) + \psi_{\text{vol}}(J) = \frac{\mu}{2} \left[I_{\widehat{C}}-3\right] + \kappa \, \mathcal{G} (J) \, .
\end{equation}
The expression for $\mathbf{S}$ is obtained from Equation \eqref{eq:Ssplit} as
\begin{equation}
    \mathbf{S} = 2 \left[\left[ \frac{\mu}{2} \, J^{-2/3} \right]\mathbf{I} + J^2 \, \left[ \kappa \left[ \frac{\partial \, \mathcal{G}}{\partial J^2} \right] - \frac{\mu}{6} \, I_{C} \, J^{-8/3} \right] \, \mathbf{C}^{-1}\right] .
\end{equation}
Choosing, $\mathcal{G} = \left[J^2-1-2\ln J\right] / 4 \, $ from Equation \eqref{eq: used G model}, one can further specialize the expression as 
\begin{equation}\label{eq: S for pstress}
    \mathbf{S} (\mathbf{C}) = \mu \, J^{-2/3} \, \mathbf{I} + \underbrace{\left[ \frac{\kappa}{2} \left[ J^2 - 1 \right] - \frac{\mu}{3} \, I_{C} \, J^{-2/3} \right]}_{\gamma= \gamma({\mathbf{C}})} \, \mathbf{C}^{-1} .
\end{equation}
Based on the decomposition in Equation \eqref{eq: c dec}, one obtains
\begin{equation}
    \mathbf{S} (\overline{\mathbf{C}}, C_{33}) = \mu \, J^{-2/3} (\overline{\mathbf{C}}, C_{33}) \, \mathbf{I} + \gamma (\overline{\mathbf{C}}, C_{33}) \ \mathbf{C}^{-1} (\overline{\mathbf{C}}, C_{33}).
\end{equation}
Given that 
$[\mathbf{C}^{-1}]_{33}=1/C_{33}$, an application of the plane stress assumption $S_{33}=0$ results in
\begin{equation}\label{eq: s33=0}
    S_{33} (\overline{\mathbf{C}}, \acute{C}_{33}) = \mu \, J^{-2/3} (\overline{\mathbf{C}}, \acute{C}_{33}) + \gamma (\overline{\mathbf{C}}, \acute{C}_{33}) \, \frac{1}{\acute{C}_{33}} = 0 ,
\end{equation}
so that by substituting $\gamma = - \mu \, \acute{C}_{33} \, J^{-2/3} (\overline{\mathbf{C}}, \acute{C}_{33})$ into Equation \eqref{eq: S for pstress}, an expression for the Piola--Kirchhoff stress specialised to plane stress is derived as
\begin{equation}
    \mathbf{S} (\overline{\mathbf{C}}, \acute{C}_{33}) = \mu \, J^{-2/3} (\overline{\mathbf{C}}, \acute{C}_{33}) \left[ \mathbf{I} - \acute{C}_{33} \, \mathbf{C}^{-1} \right] .
\end{equation}
The in-plane stress $\overline{\mathbf{S}}$ follows as
\begin{equation}\label{eq: S for pstress 2}
    \overline{\mathbf{S}} (\overline{\mathbf{C}}, \acute{C}_{33}) = \mu \, J^{-2/3} (\overline{\mathbf{C}}, \acute{C}_{33}) \left[ \overline{\mathbf{I}} - \acute{C}_{33} \, \overline{\mathbf{C}}^{-1} \right] .
\end{equation}
Using the  expression for $\overline{\mathbf{S}}$ in Equation \eqref{eq: S for pstress 2}, the first two expressions  in Equation \eqref{eq: C_IJKL for pstress} are calculated as
\begin{subequations}
\begin{align}
    \frac{\partial \overline{\mathbf{S}}}{\partial \overline{\mathbf{C}}} \bigg|_{\acute{C}_{33}} &= \acute{C}_{33} \, \mu \, J^{-2/3} \, \overline{\mathbf{C}}^{-1} \odot \overline{\mathbf{C}}^{-1} - \frac{2}{3} \, \mu \, J^{-5/3} \left[ \overline{\mathbf{I}} - \acute{C}_{33} \, \overline{\mathbf{C}}^{-1} \right] \otimes \frac{\partial J}{\partial \overline{\mathbf{C}}} \, , \\ 
    \frac{\partial \overline{\mathbf{S}}}{\partial \acute{C}_{33}} \bigg|_{\overline{\mathbf{C}}} &= - \mu \, J^{-2/3} \, \overline{\mathbf{C}}^{-1} - \frac{2}{3} \, \mu \, J^{-5/3} \left[ \overline{\mathbf{I}} - \acute{C}_{33} \, \overline{\mathbf{C}}^{-1} \right] \frac{\partial J}{\partial \acute{C}_{33}} \, ,
\end{align}
\end{subequations}
where using $J^2 = \acute{C}_{33} \, \overline{J}^2$ results in
\begin{align}
    \frac{\partial J}{\partial \overline{\mathbf{C}}} = \frac{1}{2} \, J^{-1} \, \acute{C}_{33} \, \overline{\mathbf{C}}^{-1} \, \overline{J}^2
    && \text{and} &&
    \frac{\partial J}{\partial \acute{C}_{33}} = \frac{1}{2} \, J^{-1} \, \overline{J}^2 \, .
\end{align}
The third expression in Equation \eqref{eq: C_IJKL for pstress}, $\partial \acute{C}_{33} / \partial \overline{\mathbf{C}}$, can be calculated by ensuring that $S_{33} (\overline{\mathbf{C}}, \acute{C}_{33})$ remains zero for all possible displacement changes under consideration, that is
\begin{equation}
    \frac{\text{d} {S}_{33}}{\text{d} \overline{\mathbf{C}}} = \mu \,  \frac{\text{d} J^{-2/3}}{\text{d} \overline{\mathbf{C}}} + \acute{C}^{-1}_{33} \, \frac{\text{d} \, \gamma}{\text{d} \overline{\mathbf{C}}} + \gamma \, \frac{\text{d} \acute{C}^{-1}_{33}}{\text{d} \overline{\mathbf{C}}} = \overline{\mathbf{0}} .
\end{equation}
Upon applying the chain rule to each of the terms above, one obtains
\begin{equation}
    \frac{\partial \acute{C}_{33}}{\partial \overline{\mathbf{C}}} = - \beta \left[ \frac{\partial S_{33}}{\partial \acute{C}_{33}}\bigg|_{ \overline{\mathbf{C}}} \,  \right]^{-1} ,
\end{equation}
where
\begin{align}
    \beta &= \frac{-2 \mu}{3} \, J^{-5/3} \ \overline{\mathbf{C}}^{-1} \, \overline{J}^2 + \frac{\kappa}{2} \, \overline{\mathbf{C}}^{-1} \, \overline{J}^2 - \frac{\mu}{3} \acute{C}_{33}^{-1} \left[ J^{-2/3} \, \overline{\mathbf{I}} + I_C \, \frac{\partial J^{-2/3}}{\partial \overline{\mathbf{C}}} \right] .
\end{align}
The various terms can be calculated using the following relations:
\begin{subequations}
    \begin{align}
        \frac{\partial J^{-2/3}}{\partial \overline{\mathbf{C}}} \bigg|_{\acute{C}_{33}} &= \frac{-2}{3} \, J^{-5/3} \ \overline{\mathbf{C}}^{-1} \, \overline{J}^2, \\
    \frac{\partial \gamma}{\partial \acute{C}_{33}} &= \frac{\kappa}{2} \, \overline{J}^2 - \frac{\mu}{3} \left[ I_C \, \frac{\partial J^{-2/3}}{\partial \acute{C}_{33}} + J^{-2/3} \right], \\
    \frac{\partial S_{33}}{\partial \acute{C}_{33}}\bigg|_{ \overline{\mathbf{C}}} &= \mu \, \frac{\partial J^{-2/3}}{\partial \acute{C}_{33}} + \acute{C}_{33}^{-1} \frac{\partial \gamma}{\partial \acute{C}_{33}} - \gamma \, \acute{C}_{33}^{-2} \, .
    \end{align}
\end{subequations}

By evaluating $\acute{C}_{33}$ by solving the nonlinear Equation \eqref{eq: s33=0} using a Newton-Raphson approach and updating the Piola--Kirchhoff stress and incremental constitutive tensors from Equations \eqref{eq: S for pstress 2} and \eqref{eq: C_IJKL for pstress}, a robust approach is established for implementing the plane stress condition in a finite element setting for the compressible and nearly incompressible cases.
This approach can be easily extended to accommodate other forms of energy density functions.

\subsubsection{Specialisation to the alternative neo-Hookean model}

In this section, a neo-Hookean material model that depends on the first invariant ($I_{C}$) of the total right Cauchy--Green tensor $\mathbf{C}$ is briefly discussed.
This model requires a simpler mathematical derivation than the one presented in Section \ref{sec: neo-Hookean derivations} that depends on $I_{\widehat{C}}$.
The strain energy density function is given by
\begin{equation}\label{eq:neo-hook function2}
    \psi = \frac{\mu}{2} \left[I_{C}-3-2\ln J\right] + \frac{\kappa}{2} \left[J-1\right]^2 \, .
\end{equation}
The expression for $\mathbf{S}$ is
\begin{equation}\label{eq: S for pstress2}
    \mathbf{S} (\mathbf{C}) = 2 \, \frac{\partial\psi\left(\mathbf{C}\right)}{\partial\mathbf{C}}= \mu \left[ \mathbf{I} - \mathbf{C}^{-1} \right] + \underbrace{\kappa \, J \left[ J - 1 \right]}_{\alpha=\alpha (\mathbf{C})} \, \mathbf{C}^{-1} .
\end{equation}
Similar to the derivations in Section \ref{sec: neo-Hookean derivations}, the plane stress assumption ($S_{33}=0$) leads to a calculation of $\acute{C}_{33}$.
Consequently, the expression for the in-plane Piola--Kirchhoff stress 
specialised to plane stress is
\begin{equation}\label{eq: S for pstress 2b}
    \overline{\mathbf{S}} 
    = \mu \left[ \overline{\mathbf{I}} - \acute{C}_{33} \, \overline{\mathbf{C}}^{-1} \right] .
\end{equation}
Thus, the first two expressions  in Equation \eqref{eq: C_IJKL for pstress} are calculated as
\begin{align}
    \frac{\partial \overline{\mathbf{S}}}{\partial \overline{\mathbf{C}}} \bigg|_{\acute{C}_{33}} = \acute{C}_{33} \, \mu \, \overline{\mathbf{C}}^{-1} \odot \overline{\mathbf{C}}^{-1} && \text{and} && 
    \frac{\partial \overline{\mathbf{S}}}{\partial \acute{C}_{33}} \bigg|_{\overline{\mathbf{C}}} = - \mu \, \overline{\mathbf{C}}^{-1} \, .
\end{align}
Finally, the third term in Equation \eqref{eq: C_IJKL for pstress}, $\partial \acute{C}_{33} / \partial \overline{\mathbf{C}} \,$, is derived by enforcing the condition $S_{33} (\overline{\mathbf{C}}, \acute{C}_{33})$ remains zero for all displacement variations.
This leads to
\begin{equation}
    \frac{\text{d} {S}_{33}}{\text{d} \overline{\mathbf{C}}} = \acute{C}^{-1}_{33} \, \frac{\text{d} \alpha}{\text{d} \overline{\mathbf{C}}} + [\alpha-\mu] \, \frac{\text{d} \acute{C}^{-1}_{33}}{\text{d} \overline{\mathbf{C}}} = \overline{\mathbf{0}} .
\end{equation}
Applying the chain rule and simplifying, one obtains
\begin{equation}
    \frac{\partial \acute{C}_{33}}{\partial \overline{\mathbf{C}}} = - \acute{C}_{33}^{-1} \, \frac{\partial \alpha}{\partial \overline{\mathbf{C}}} \left[ \frac{\partial S_{33}}{\partial \acute{C}_{33}}\bigg|_{ \overline{\mathbf{C}}} \,  \right]^{-1} ,
\end{equation}
where 
\begin{equation}
    \frac{\partial \alpha}{\partial \overline{\mathbf{C}}} = \kappa \, \acute{C}_{33} \, \overline{\mathbf{C}}^{-1} \, \det \overline{\mathbf{C}} \, \left[1-\frac{1}{2} J^{-1}\right] .
\end{equation}

This alternative neo-Hookean model clearly offers a more straightforward formulation compared to the decoupled model presented in the previous section and a similar form has been used by  \citet{pascon2019large}.
However, the decoupled model in Section \ref{sec: neo-Hookean derivations} is preferred for modelling nearly incompressible materials \cite{holzapfel2000nonlinear}.

\subsection{Flatland (and plane strain)}

\paragraph{Flatland:}
In this approach, the real-world complexity is condensed into a truly two-dimensional representation by assuming the existence of only two dimensions.
This assumption is valid, for example, in the case of one-atom-layer materials made from carbon.
One can adapt three-dimensional formulations to two dimensions, by simply reducing the dimensions of the mathematical model.
The components of the deformation gradient in flatland, $\overline{\mathbf{F}} \in \mathbb{R}^{2} \otimes \mathbb{R}^{2}$, are given as
\begin{equation}
    [\overline{\mathbf{F}}] = \left[\begin{matrix} F_{11}&F_{12} \\ 
            F_{21}&F_{22} \\
    \end{matrix}\right] .
\end{equation}
The determinant of the flatland deformation gradient is denoted by $\overline{J} = \det{\overline{\mathbf{F}}}$.
For the examples in Section~\ref{sec: numerical}, a two-dimensional neo-Hookean model for hyperelastic modelling is employed where
\begin{equation}
    \psi = \frac{\mu}{2} \left[{I \,}_{\widehat{\overline{C}}}-2\right] + \kappa \, \mathcal{G} ( \, \overline{J} \, ) \, .
\end{equation}

\paragraph{Plane strain:} 
For completeness and comparison we briefly outline the plane strain case in which the deformation along the $X_3$ coordinate is constrained to zero, thus the components of the deformation gradient $\mathbf{F}$ and consequently the Green--Lagrange tensor $\mathbf{E}$ are given as
\begin{align}
    [\mathbf{F}]=\left[\begin{matrix} F_{11}&F_{12}&0 \\ 
                F_{21}&F_{22}&0 \\
                0&0&1 \\
    \end{matrix}\right]
    && \text{and} &&
    [\mathbf{E}]= \left[\begin{matrix} E_{11}&E_{12}&0 \\ 
     E_{21}&E_{22}&0 \\
     0&0&0 \\
     \end{matrix}\right].
\end{align}

\section{Numerical results}\label{sec: numerical}

The finite element formulation developed here, implemented using the deal.II library, version 9.2 \citep{arndt2020deal}, is available at \citep{madeal2024}.
Deal.II is a versatile C++ library for the finite element method, facilitating the efficient development of modern finite element codes through its object-oriented architecture \citep{bangerth2007deal, bangerth2007deal2}.

Computations for various nonlinear problems with different loading scenarios are conducted to investigate the performance of the different approaches with an emphasis on the plane stress approximation.
Two well-known benchmark problems from the literature, namely, Cook's cantilever and an inhomogeneous compression problem, are examined to verify the accuracy of the models, comparing the results with prior studies.
Further examples, including stretching of composites reinforced by particles and fibres, are presented to compare the performance of the models in evaluating the extreme nonlinear behaviour of these multi-phase structures.
For all examples, the neo-Hookean material model with the decoupled free energy function \eqref{eq:neo-hook function} and volumetric contribution \eqref{eq: used G model} 
are selected.
The number of equal load increments is set to $10$ for all problems. 
Quadrilateral elements with polynomial order $p_o$ are used for approximating the finite element basis functions in two dimensions, while hexahedral elements 
are employed for three-dimensional problems.
The stress values depicted in the contour plots are averaged element stresses.
The effective von--Mises stress $\sigma_{\text{eff}}$ is given by
\begin{equation}
    \sigma_{\text{eff}} = \sqrt{\frac{3}{2} \, \bm{\sigma}_{\text{dev}} : \bm{\sigma}_{\text{dev}} } \ , \quad \text{where,} \quad \bm{\sigma}_{\text{dev}} = \bm{\sigma} - \frac{1}{d} \, \text{tr} (\bm{\sigma}) \, \mathbf{I} ,
\end{equation}
where $d=3$ for three dimensions and plane stress/strain and $d=2$ for the flatland approach.
Since it is more convenient to represent material incompressibility using the Poisson's ratio $\nu$, the material properties are expressed in terms of $\nu$ and the shear modulus $\mu$, where the Poisson's ratio is given by $\nu=[{3 \kappa - 2 \mu}]/\big[{2[3 \kappa + \mu]}\big]$.

The consistent linearization of the plane stress implementation results in quadratic convergence of both the inner and outer Newton schemes.

\subsection{Cook's cantilever}

\begin{figure}
    \centering
    \begin{tabular}{m{0.45\textwidth} m{0.45\textwidth}}
        \includegraphics[width=\linewidth]{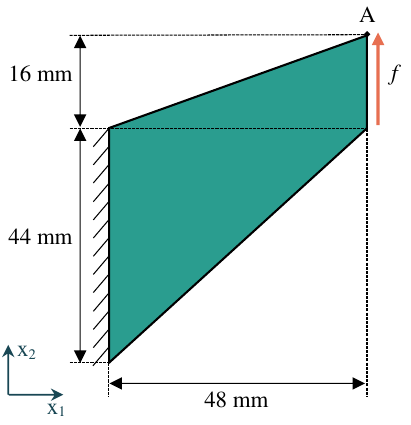} & \input{ai.tex} \\
        \centering (a) & \centering (b)
    \end{tabular}
    \caption{The Cook's cantilever under transverse loading: (a) Schematic representation. (b) Vertical displacement of the tip versus applying traction $f$.}
    \label{fig:cook}
\end{figure}

The first example is the 
bending of the Cook's cantilever \citep{cook1977ways}.
This problem is associated with a significant amount of shear deformation.
A comprehensive analysis is presented, and the results compared with those found in prior studies.
Figure \ref{fig:cook} (a) depicts the geometry of a Cook's cantilever under a traction $f$ at the right boundary.
The thickness of the cantilever is $t = 1 \, \text{mm}$.

The material properties
are a shear modulus of $\mu = 80.1938 \, \text{MPa}$, and Poisson's ratio of $\nu = 0.4999 \,$, chosen in accordance with \citep{simo1992geometrically, betsch19964, wriggers1996formulation, reese2002equivalent, sze2004stabilized, angoshtari2017compatible, pascon2019large, dhas2022mixed}. 
The geometry is uniformly meshed using a grid of size $n \times n$, where $n$ takes values from the set $\{2, 4, 8, 16, 32, 64\}$ representing the number of elements per edge.
The three-dimensional models are meshed with the same mesh size with one element through the thickness, and the mesh size denoted as $n \times n \times 1 \,$.
The three-dimensional model is unconstrained in the $x_3$ direction, allowing free movement along this axis.

The vertical displacements $u_2$ of the cantilever tip (point A in Figure \ref{fig:cook}a) calculated for the traction values of $f=24 \, \text{N}/\text{mm}^2$ and $f=40 \, \text{N}/\text{mm}^2$ are given in Table \ref{tab:cook1}.
This table compares different finite element models with linear ($p_o=1$) and quadratic ($p_o=2$) basis functions, as well as various mesh sizes.

\begin{table}
    \centering
    \caption{Vertical displacement $u_2$ (mm) of the point A of Cook's cantilever.}
    \begin{tabular}{|c|c|c|c|c|c|c|c|c|c|}
         \hline
         Applied $f$ & \multirow{2}{*}{Mesh} & \multicolumn{2}{c|}{flatland} & \multicolumn{2}{c|}{plane strain} & \multicolumn{2}{c|}{plane stress} & \multicolumn{2}{c|}{three dim.*} \\
         \cline{3-10}
         ($\text{N}/\text{mm}^2$) & & $p_{o}=1$ & $p_{o}=2$ & $p_{o}=1$ & $p_{o}=2$ & $p_{o}=1$ & $p_{o}=2$ & $p_{o}=1$ & $p_{o}=2$ \\
         \hline
         $24$ & $2 \times 2$ & 13.91 & 18.45 & 13.77 & 18.29 & 14.42 & 19.26 & 15.47 & 19.85 \\
          & $4 \times 4$ & 16.69 & 18.20 & 16.65 & 18.17 & 17.85 & 19.75 & 18.40 & 19.89 \\
          & $8 \times 8$ & 17.70 & 18.20 & 17.68 & 18.18 & 19.24 & 19.88 & 19.47 & 19.93 \\
          & $16 \times 16$ & 18.00 & 18.20 & 17.99 & 18.19 & 19.71 & 19.93 & 19.79 & 19.95 \\
          & $32 \times 32$ & 18.11 & 18.21 & 18.10 & 18.20 & 19.87 & 19.95 & 19.89 & 19.96 \\
          & $64 \times 64$ & 18.15 & 18.22 & 18.15 & 18.21 & 19.93 & 19.97 & 19.93 & 19.96 \\
        \hline
        $40$ & $2 \times 2$ & 19.56 & 24.27 & 19.67 & 24.16 & 21.23 & 25.17 & 21.73 & 25.78 \\
         & $4 \times 4$ & 22.40 & 24.20 & 22.46 & 24.17 & 24.18 & 25.78 & 24.42 & 26.01 \\
         & $8 \times 8$ & 23.54 & 24.23 & 23.55 & 24.22 & 25.29 & 25.98 & 25.50 & 26.10 \\
         & $16 \times 16$ & 23.93 & 24.25 & 23.93 & 24.25 & 25.75 & 26.07 & 25.88 & 26.15 \\
         & $32 \times 32$ & 24.08 & 24.28 & 24.10 & 24.27 & 25.96 & 26.13 & 26.02 & 26.17 \\
         & $64 \times 64$ & 24.16 & 24.31 & 24.16 & 24.30 & 26.06 & 26.16 & 26.09 & 26.18 \\
        \hline
        \multicolumn{10}{l}{*The mesh for three-dimensional models incorporates one element thickness in the $x_3$ direction}
    \end{tabular}
    \label{tab:cook1}
\end{table}

The results indicate that the flatland and plane strain models are nearly identical when the mesh is sufficiently refined.
The plane strain case is verified using \citep{reese2002equivalent, angoshtari2017compatible, dhas2022mixed}:
the vertical displacement of the tip of the Cook's cantilever for $f=24 \, \text{N}/\text{mm}^2$ is reported as $u_2=18.05 \, \text{mm}$ by \citet{reese2002equivalent}, $u_2=18.2 \, \text{mm}$ by \citet{angoshtari2017compatible} and $u_2=18.1 \, \text{mm}$ by \citet{dhas2022mixed} for the finest mesh used in their works.
\citet{pascon2019large} uses a plane stress approximation with identical material properties and the traction of $f=40 \, \text{N}/\text{mm}^2$.
They obtained displacements at the tip of $u_1 = -28.12 \, \text{mm}$ and $u_2 = 26.22 \, \text{mm}$ for their finest mesh.
In the present work, the finest mesh of $64 \times 64$ and $p_o=2$ yields displacements of the tip of $u_1 = -28.04 \, \text{mm}$ and $u_2 = 26.16 \, \text{mm}$, showing good agreement with a difference of less than $0.3 \%$ compared to their results.


It is evident that the results from the plane stress and three-dimensional models align closely.
Vertical displacements for the plane stress and three-dimensional cases are notably higher than those for the flatland and plane strain models.
Furthermore, it is apparent that results for higher-order elements converge more rapidly, with this effect being more pronounced than that of mesh refinement.

It should be noted that the three-dimensional model is unconstrained in the $x_3$ direction.
However, if the displacement of the cantilever is fixed in the $x_3$ direction so that it is trapped between two walls, the results align more closely with the flatland and plane strain cases.
In this case, using the most refined mesh ($64 \times 64 \times 1$) with $p_o=2$, the vertical displacement of the tip for the three-dimensional model when applying $f=24 \, \text{N}/\text{mm}^2$ is calculated as $u_2 = 18.21 \, \text{mm}$, which is closer to the plane strain case. 

To further verify the implementation, the vertical displacement of the cantilever tip versus the applied traction $f$ is depicted in Figure \ref{fig:cook} (b).
It should be noted that the results from other sources correspond to the finest mesh used in those studies.
This graph highlights the nonlinear relationship between the applied load and the displacement of the tip.
%
Note that some of the references that the present results were compared against used a slightly different version of the volumetric energy density function, such as the ones shown in Equations \eqref{eq:gfunctions}.
However, despite this variation, the differences observed in the results are negligible since the material is nearly incompressible ($\nu = 0.4999$).

\begin{figure}
	\begin{center}
		\begin{tabular}{c c c c c}
            flatland & plane strain & plane stress & three dimensions \\
		    \includegraphics[scale=0.22]{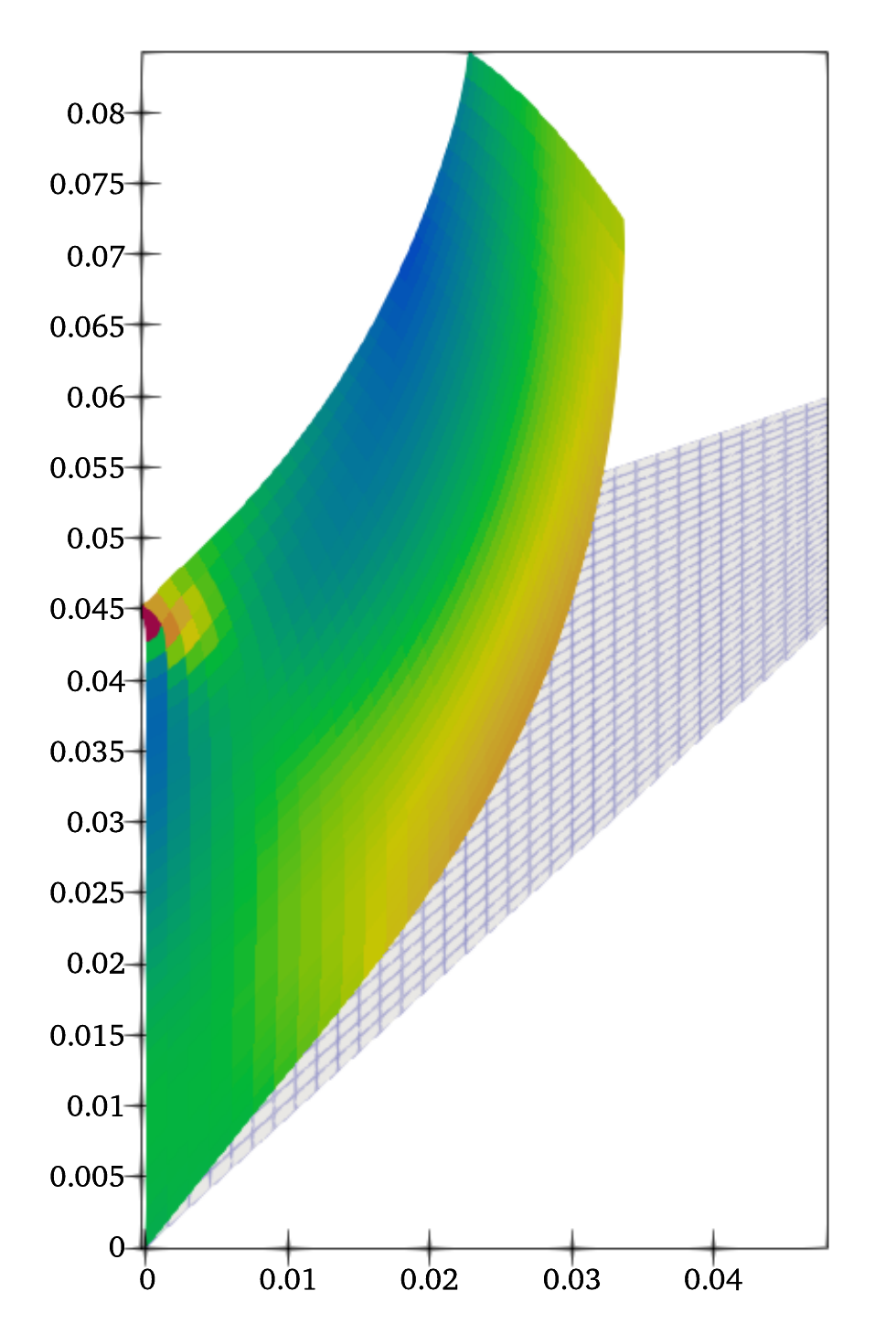} & \includegraphics[scale=0.22]{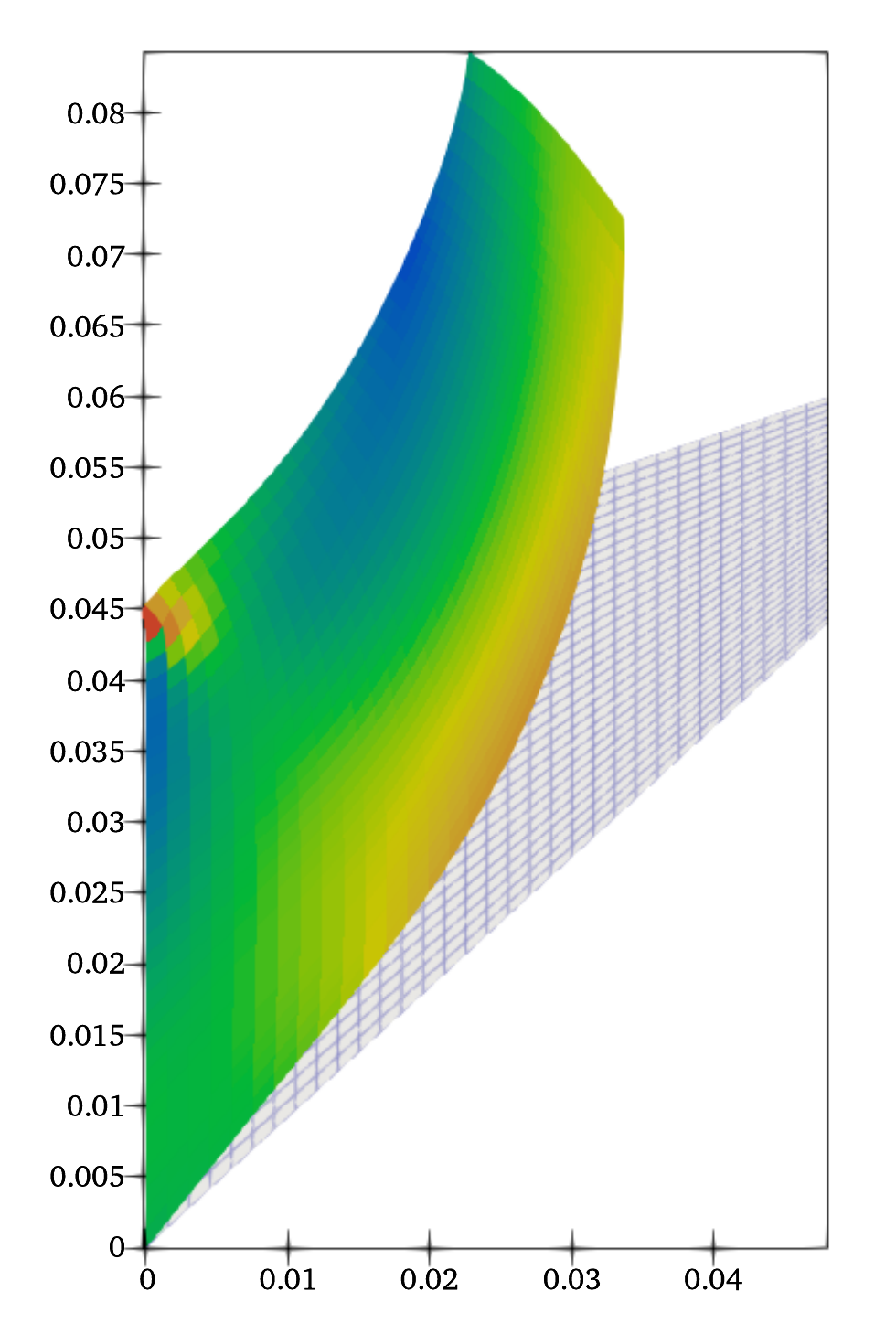} &
            \includegraphics[scale=0.22]{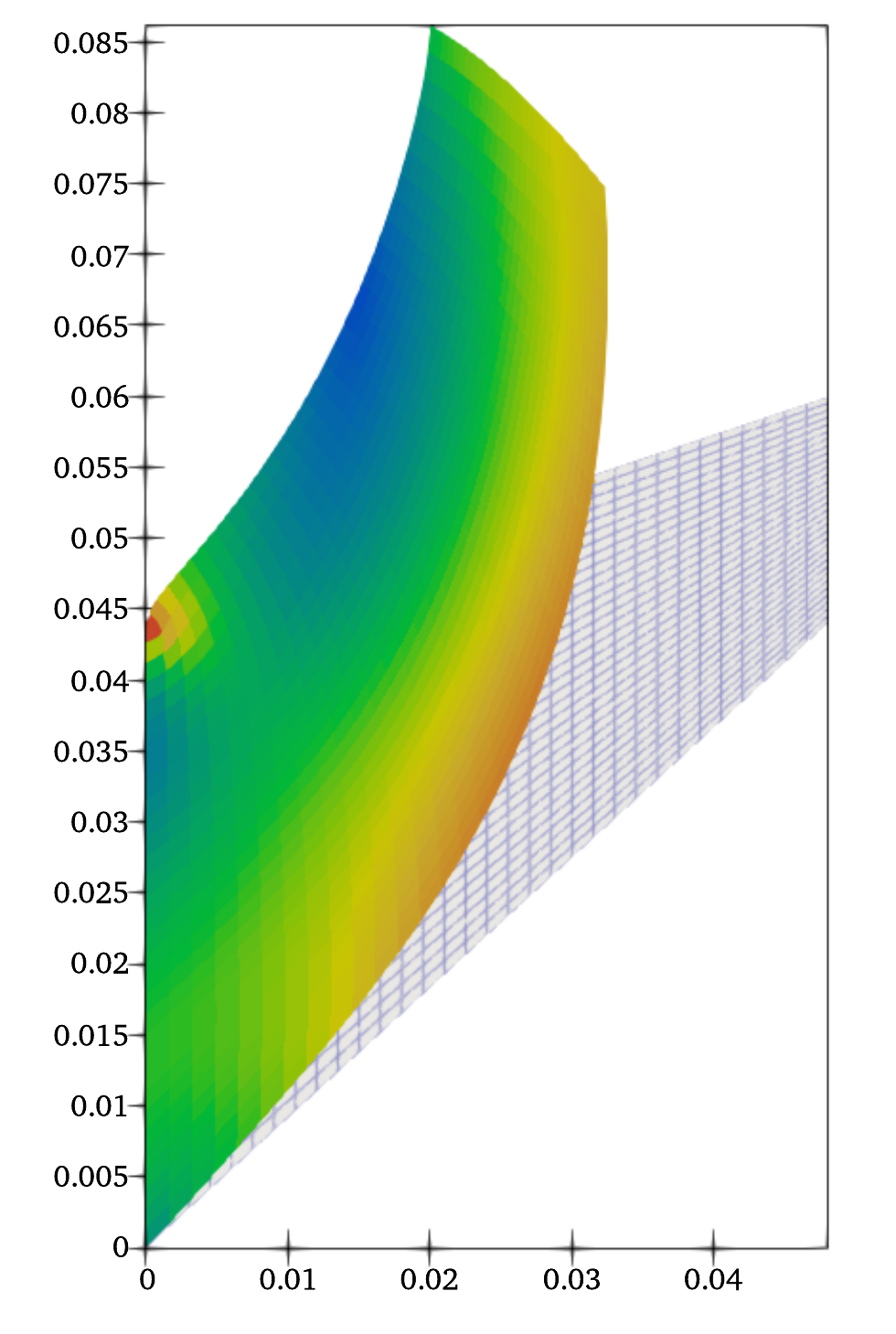} & \includegraphics[scale=0.22]{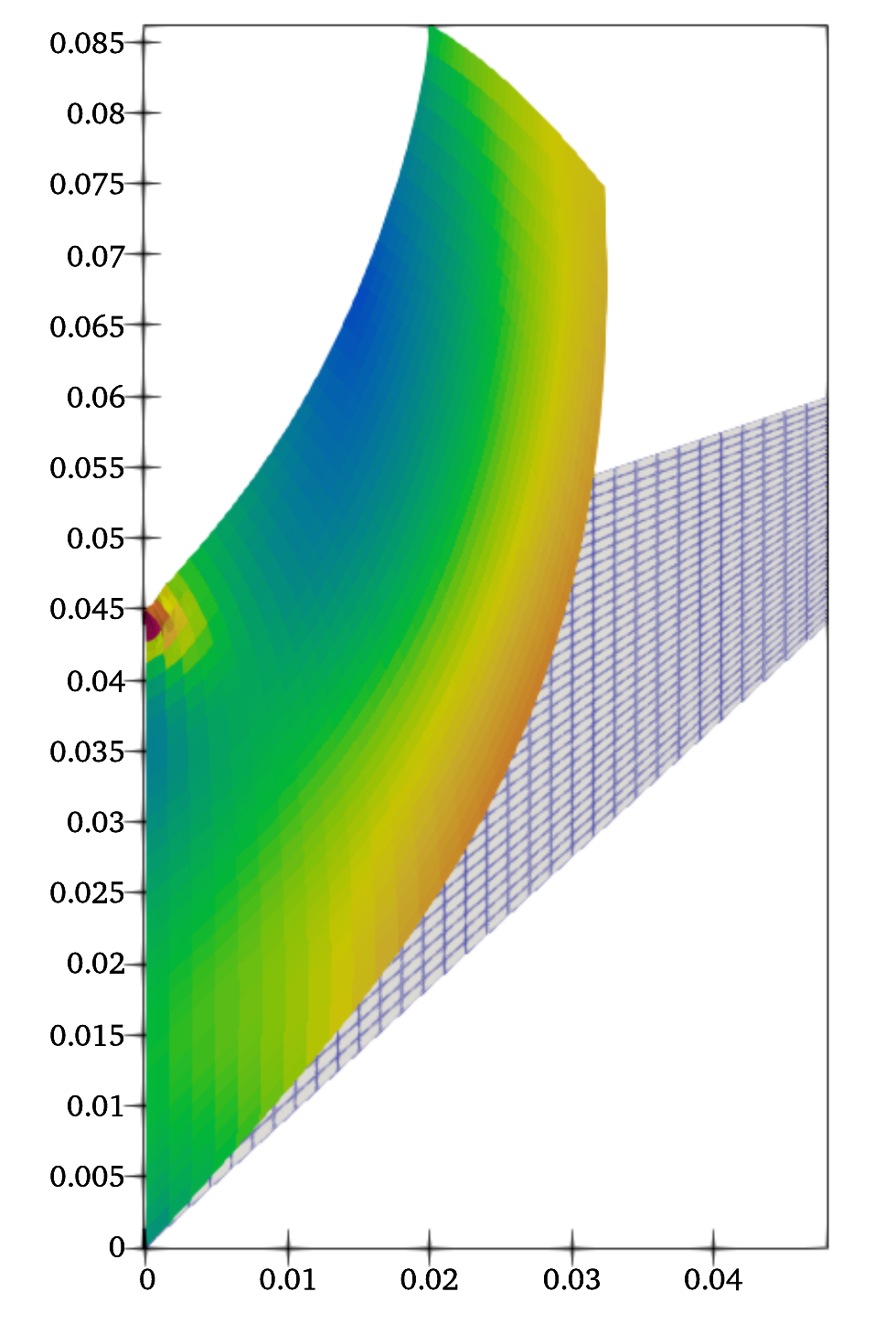} & \includegraphics[scale=0.22]{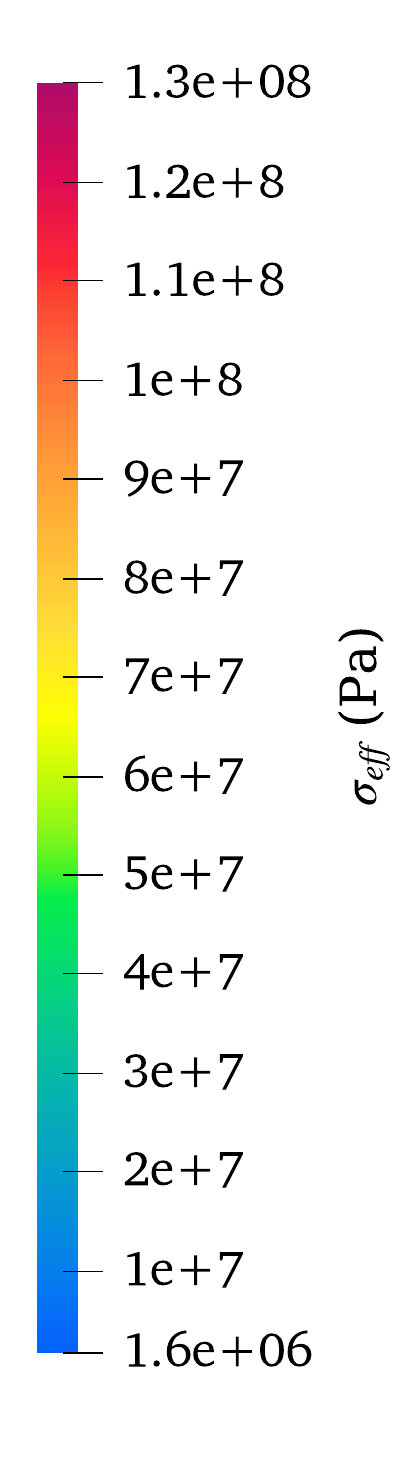}
		\end{tabular}
	\end{center}
    \caption{The deformed Cook's cantilever with $32 \times 32$ mesh size, $p_o=2$ for different models under $f=40 \, \text{N}/\text{mm}^2$.}
    \label{fig:cook contour}
\end{figure}

For the sake of illustration, Figure \ref{fig:cook contour} shows the deformed Cook's cantilever  for the $32 \times 32$ mesh, $p_o=2$ and $f=40 \, \text{N}/\text{mm}^2$ for different configurations.
The first element from the top left corner is critical in this problem since it undergoes a significant distortion.
As seen from the contours, this element is more distorted in flatland, plane strain and three-dimensional models than in the plane stress model. 

\begin{figure}
	\begin{center}
		\begin{tabular}{c}
			\input{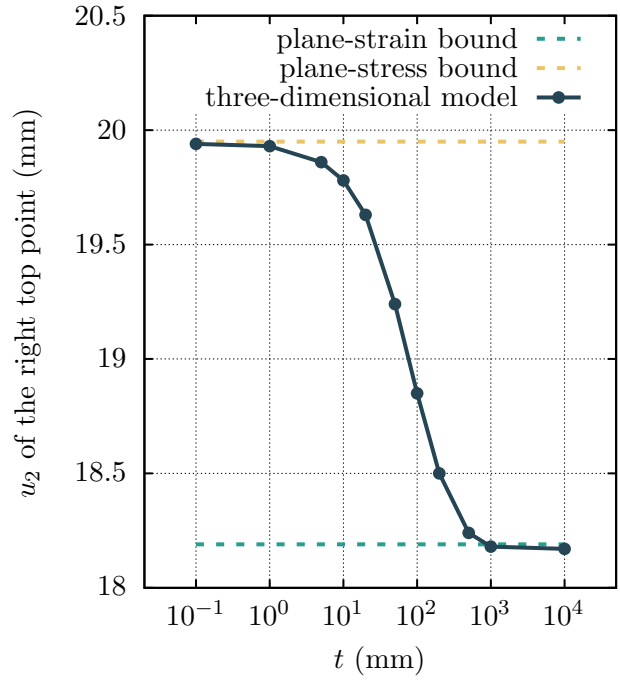}
		\end{tabular}
	\end{center}
 \caption{The effect of the thickness of the Cook's cantilever in the three-dimensional model on the vertical displacement of the tip.}
 \label{fig:cook thickness}
\end{figure}

Figure \ref{fig:cook thickness} illustrates the impact of the thickness $t$ of the cantilever in the three-dimensional model on the vertical displacement of the structure's tip for $f=[24 \, t] \, \text{N}/\text{mm}^2$. 
As the thickness increases, the results from the three-dimensional models diverge further from the plane stress bound and approach those of the flatland/plane strain models, as anticipated.
This highlights the importance of carefully selecting modelling approaches, considering the nature of the problem under analysis.

Analysing the computational cost of the different models, the following trend is observed: $t_{\text{three-dimensional}} \gg t_{\text{plane strain}} > t_{\text{flatland}} > t_{\text{plane stress}} \, $.
The number of elements for the finest mesh, $32 \times 32$, is 1024 for all four models, while the number of nodes are 14594, 14594, 8450, and 46217, and the total number of degrees of freedom are 7297, 7297, 4225, and 15405, for the flatland, plane strain, plane stress, and three-dimensional models with $p_o=2$, respectively.
As expected, the computational time for the three-field mixed-formulation in the three-dimensional model is higher than that of the plane strain model, which, in turn, is higher than that of the flatland model.
The computational time for the classic single-field model for plane stress is lower than that of the plane strain mixed-formulation model, despite having additional Newton iteration required to determine $C_{33} \,$.
As expected, the plane stress model exhibits significantly lower runtime compared to the three-dimensional model.
For example, for the finest mesh with $p_o=2$, the runtime for the three-dimensional model is approximately 50 times longer than that of the plane stress model.

\subsection{Inhomogeneously compressed block}

\begin{figure}[http]
    \centering
    \includegraphics{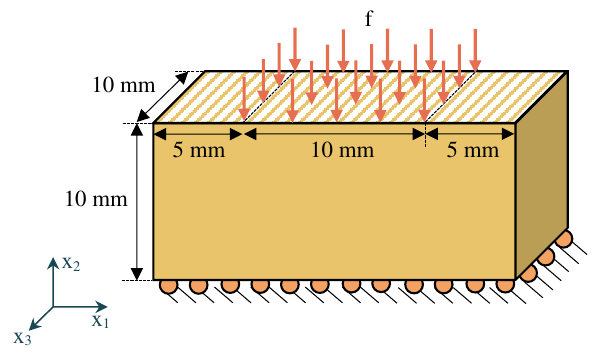}
    \caption{The block under inhomogeneous compression; the top surface is constrained in $x_1$ direction.}
    \label{fig:indentation}
\end{figure}

The second example is that of a nonlinear elastic block undergoing inhomogeneous compression, also known as the punch problem.
Figure \ref{fig:indentation} illustrates the geometry, loading, and boundary conditions.
The bottom surface is fixed in the $x_2$ direction but can freely move in the $x_1$ direction.
It is also noted that the top surface of the block is constrained in the $x_1$ direction.
For the three-dimensional model, the back and front surfaces are constrained in the $x_3$ direction, simulating a scenario where the block is confined between two rigid walls.
The material properties are the same as the previous example.
To leverage problem symmetry, the finite element analysis only considers half of the geometry by dividing the block into two halves, each measuring $10 \, \text{mm} \times 10 \, \text{mm} \times 10 \, \text{mm}$.

\begin{figure}
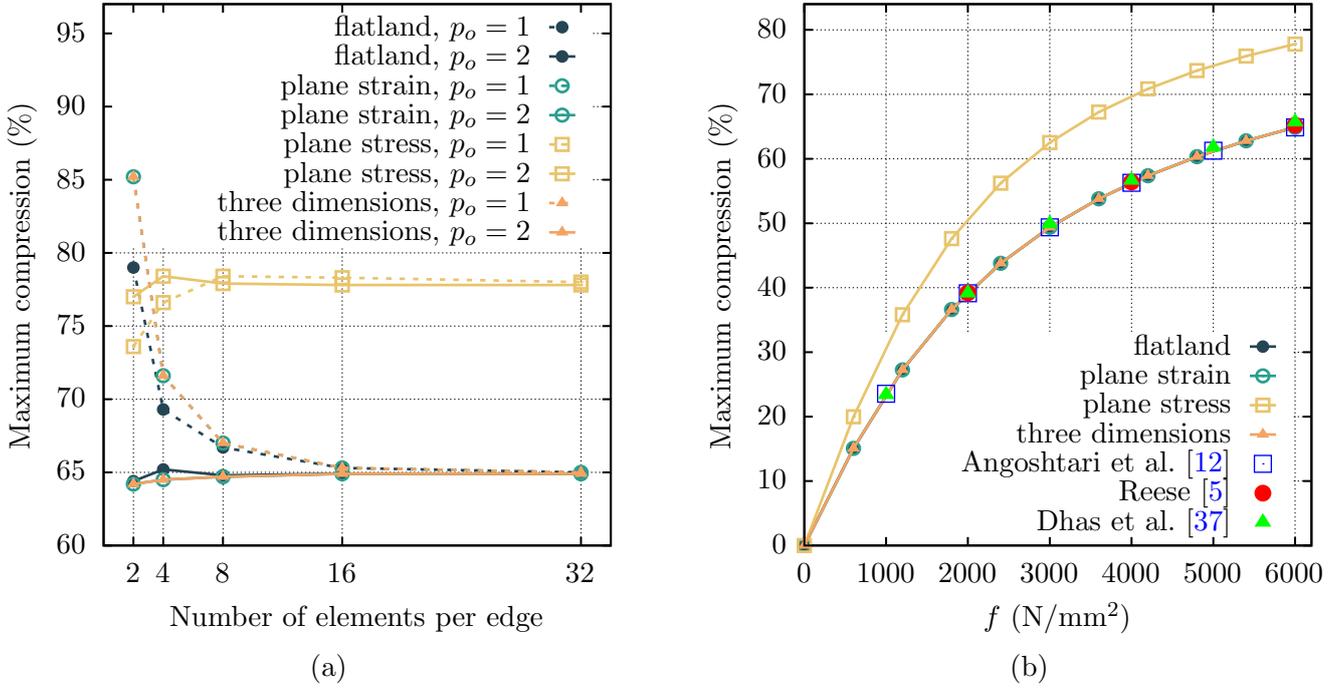

	\begin{center}
		\begin{tabular}{c c}
			\input{b1.tex} & \input{b2.tex} \\
            (a) & (b)
		\end{tabular}
	\end{center}
 \caption{Maximum compression of the block: (a) for different mesh sizes and applied traction of $f=6000 \, \text{N}/\text{mm}^2$ (b) versus applied traction $f$.}
 \label{fig:punch plot1}
\end{figure}

The geometric domain is uniformly meshed, that is, $n \times n$ ($n \times n \times n$ for three dimensions), where $n$ is the number of elements per each edge and $n\in\{2,4,8,16,32\}$.
The maximum compression of the block is defined as the percentage of the maximum vertical displacement of the middle node of the top-surface relative to the height of the block.
Figure \ref{fig:punch plot1} (a) illustrates the maximum compression across different models, considering various meshes and polynomial degrees of the elements (basis functions).
By $n=16$, the results appear to be converged across all different models. 
Moreover, models utilizing quadratic elements exhibit faster convergence than those with linear elements; even a model with $n=2$ and quadratic elements demonstrates high accuracy.
Given that the block is confined between two walls, restricting movement in the $x_3$ direction, the three-dimensional model is expected to converge to the plane strain model, which indeed occurs.
Interestingly, the block under plane stress condition experiences $\sim 20 \%$ higher compression than other models.

Figure \ref{fig:punch plot1} (b) depicts the maximum compression of the block versus the applied traction $f$.
The maximum compression demonstrates nonlinear behaviour across the range of applied traction.
The recurring trend persists; the plane stress model demonstrates higher compression relative to all other models, a distinction notably apparent even at lower loads.
The results also validate the (three-dimensional, flatland, and the) plane strain case by comparison against other studies \citep{angoshtari2017compatible, reese2002equivalent, dhas2022mixed}, demonstrating satisfying agreement.

\begin{figure}[htbp]
	\begin{center}
		\begin{tabular}{c c}
            plane strain & plane stress \\
			\includegraphics[scale=0.19]{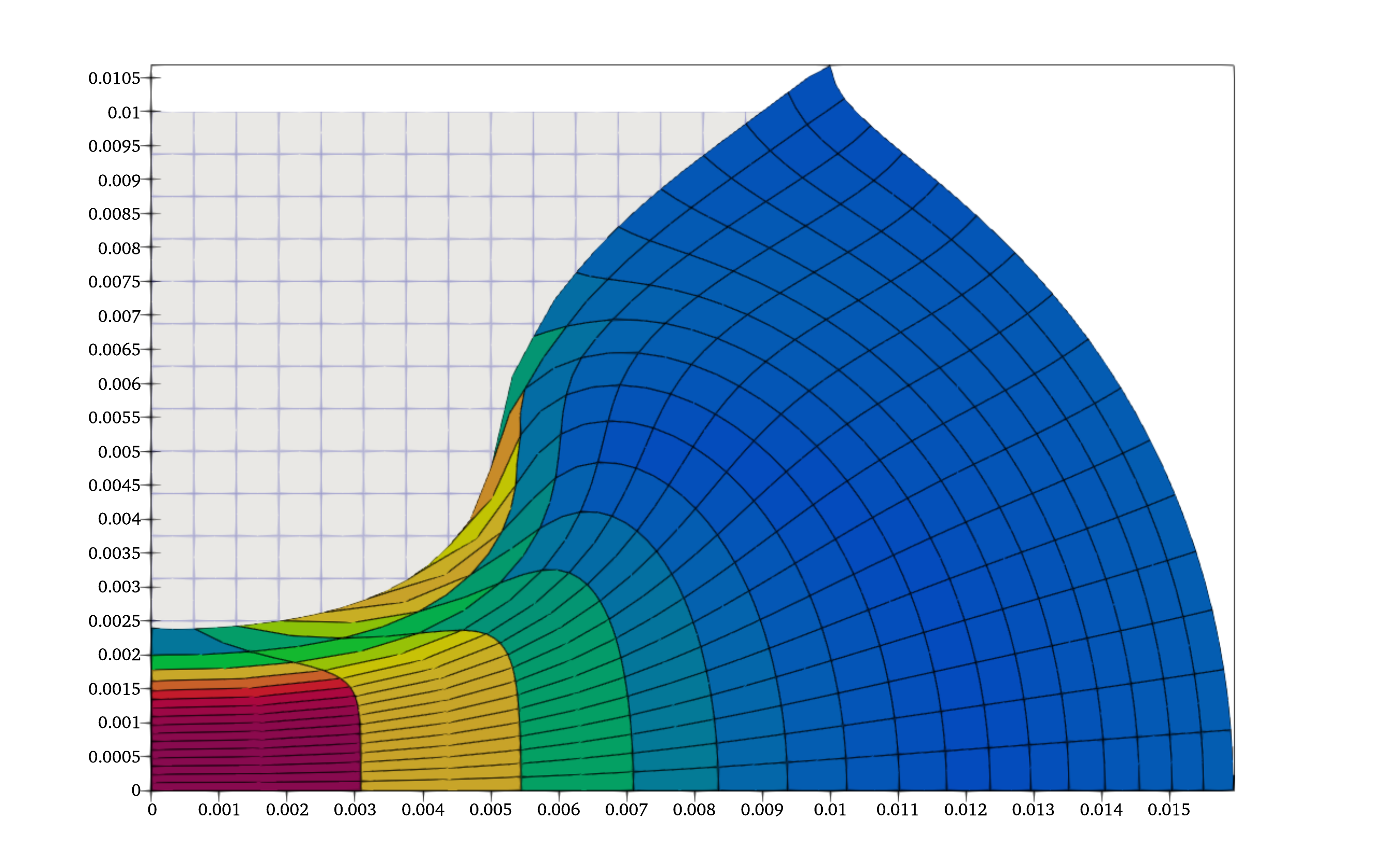} 
            & \includegraphics[scale=0.19]{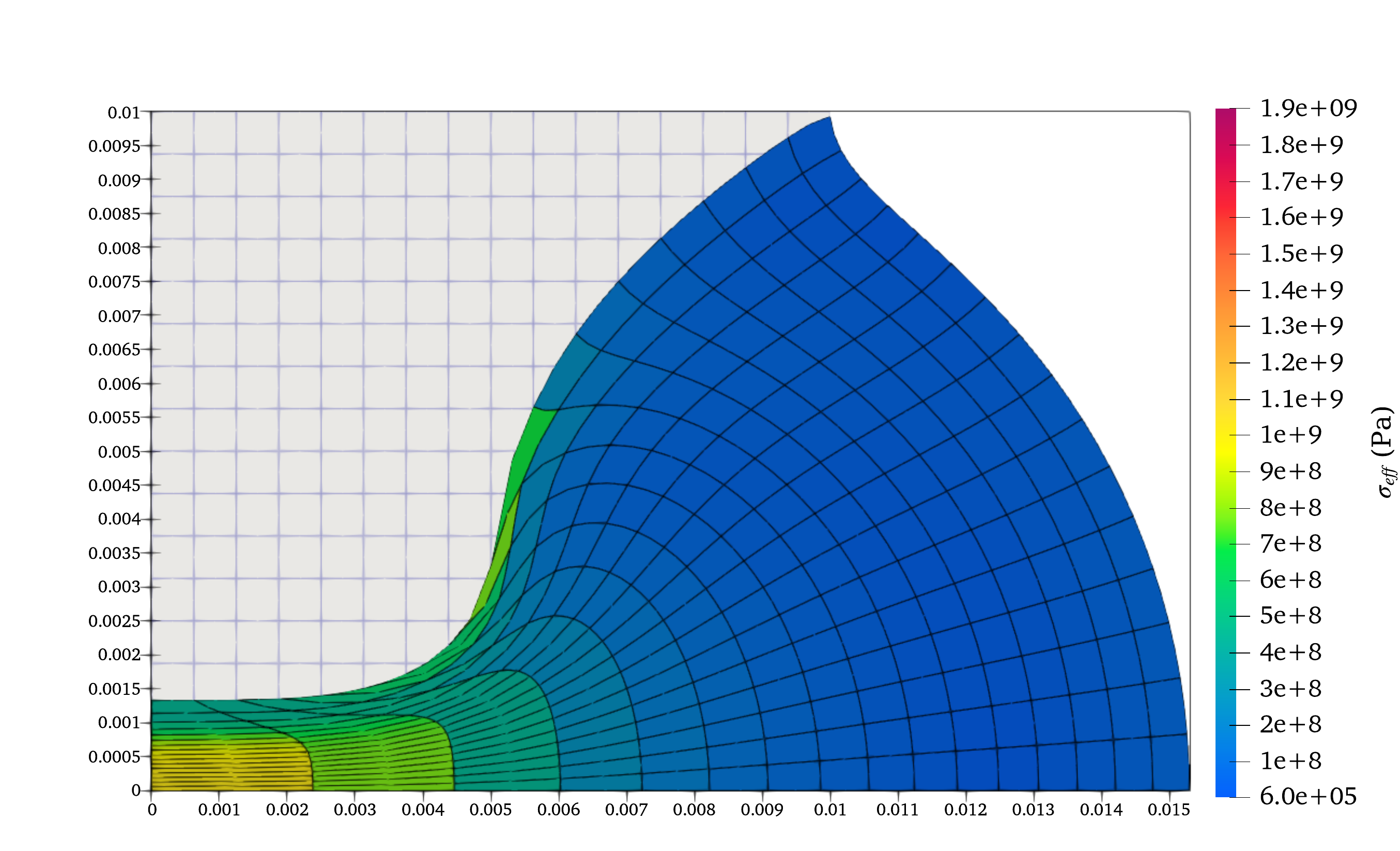}
		\end{tabular}
	\end{center}
 \caption{Compressed block with applied traction of $f=12000 \, \text{N}/\text{mm}^2$ under planar conditions.}
 \label{fig:punch compressed}
\end{figure}

To validate the results of the plane stress model, a traction of $f=12000 \, \text{N}/\text{mm}^2$ is applied, following the study by \citet{pascon2019large}.
The maximum compression under this traction is measured as $86.7 \%$, closely matching the approximately $87 \%$ compression reported by \citet{pascon2019large} at the centre. 
The compressed block with $n=16$ and $p_o=2$ 
is illustrated in Figure \ref{fig:punch compressed}, demonstrating the capability of the formulation to model extreme deformations.

\subsection{Stretching of heterogeneous composites}

In the final example, the complex deformation behaviour of two square blocks containing several circular particles and short fibres, scenarios commonly encountered in composite materials and reinforced structures \citep{ahmadi2019micromechanical, ahmadi2017multi, gao2019numerical}, is investigated.
The left end of the domain is fixed, and the right end is displaced to a length twice that of the block.
The matrix material is characterized by a shear modulus of $\mu_m = 1 \, \text{MPa}$ and a Poisson's ratio of $\nu_m=0.4999$ for the nearly incompressible case and $\nu_m=0.3$ for the compressible case.
The particles and fibres possess a shear modulus of $\mu_p=50 \, \text{MPa}$ and a Poisson's ratio of $\nu_p=0.3 \, $.
Sufficiently refined meshes are utilized to discretise the geometry, ensuring the convergence of results.

\begin{figure}[htb]
	\begin{center}
		\begin{tabular}{c c}
            plane stress compressible & plane stress nearly incompressible \\
		    \includegraphics[scale=0.21]{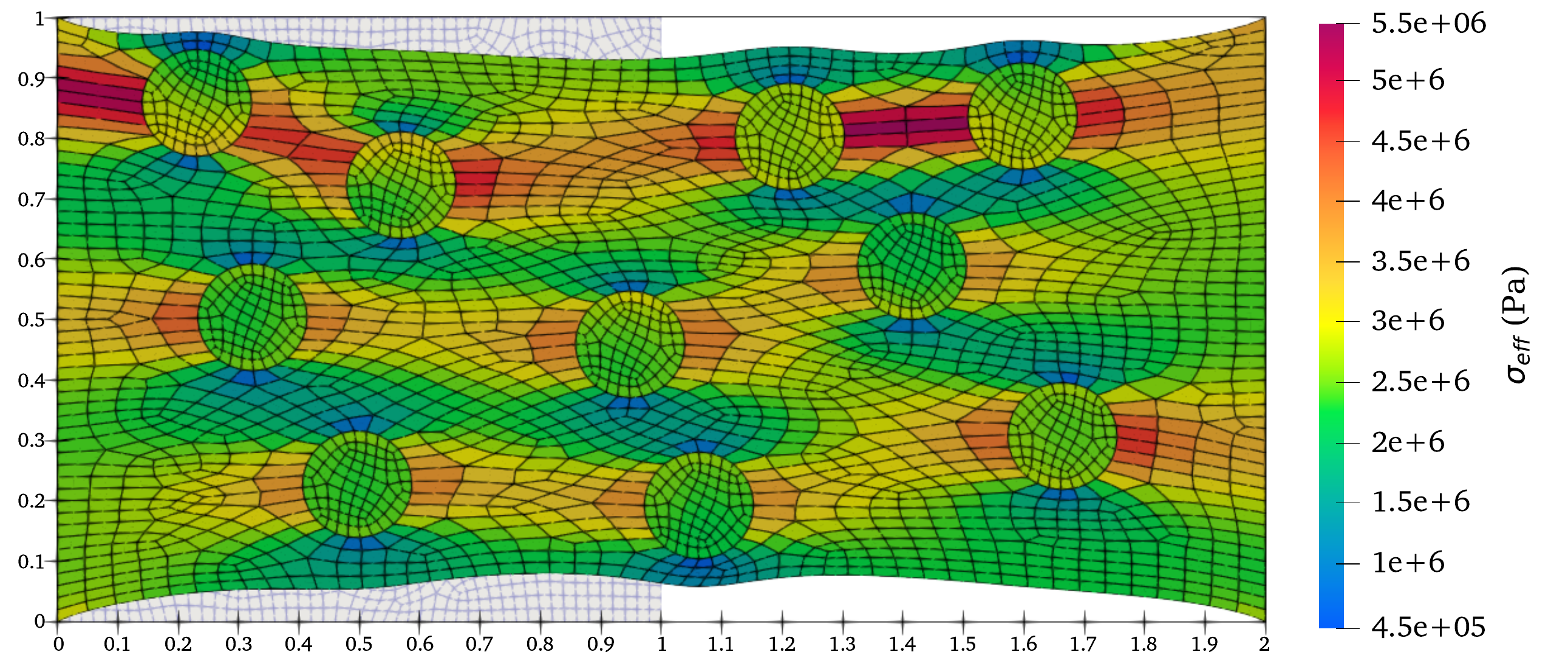} & \includegraphics[scale=0.21]{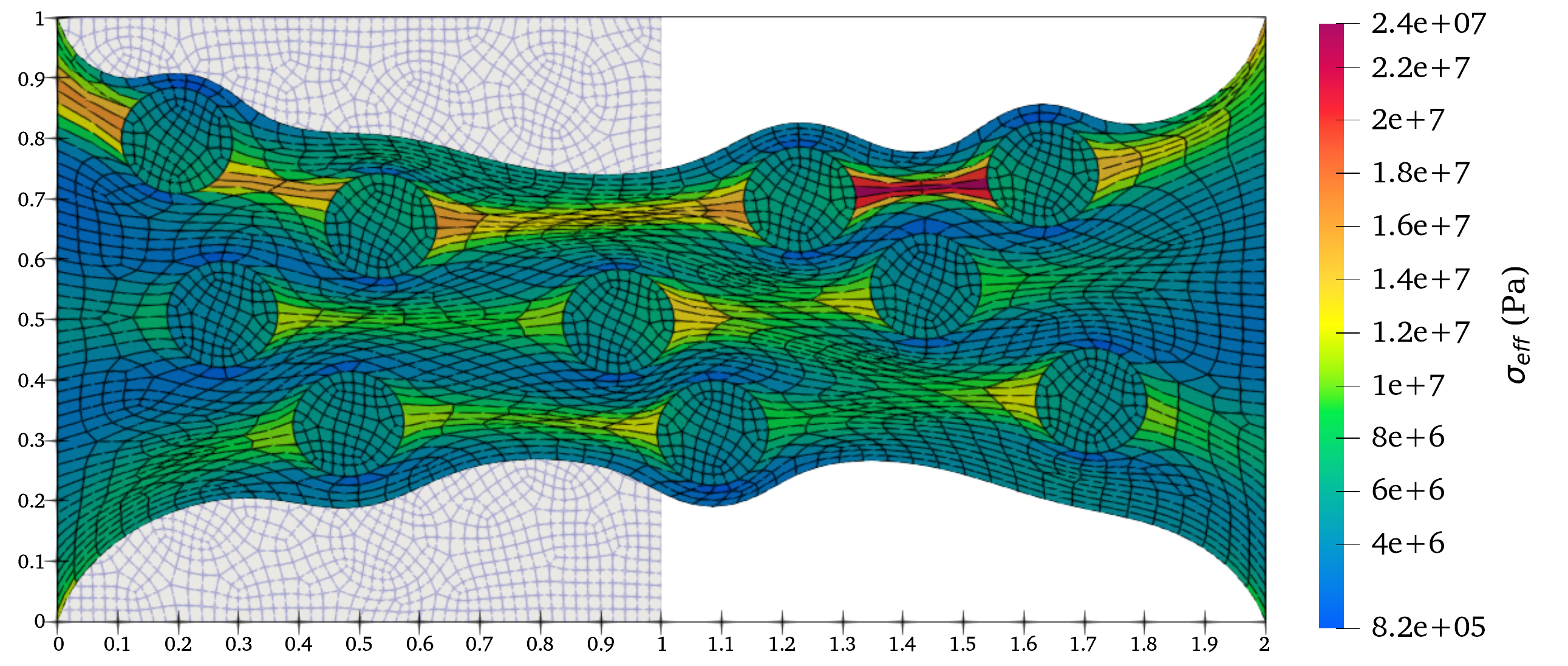}
		\end{tabular}
	\end{center}
    \caption{The plane stress compressible ($\nu=0.3$) and nearly incompressible ($\nu=0.4999$) stretched composite with several particles.}
    \label{fig:pwp_contours}
\end{figure}

In the first scenario, the volume fraction of the particles is fixed at $25 \%$ of the total volume of the composite in the material configuration.
Ten particles, with stiffness greater than the matrix, are randomly positioned within the matrix material.
Figure \ref{fig:pwp_contours} visualizes the deformed composites reinforced by circular particles under plane stress conditions for compressible and nearly incompressible materials.
The contour bars illustrate the distribution of von--Mises stress values.
The results indicate that compressible composite exhibits less deformation perpendicular to the stretch direction compared to the nearly incompressible one. 
This behaviour arises due to the compressible nature of the material, which allows for more volumetric changes. 
Additionally, the circular particles mostly remain undeformed due to their higher relative stiffness. 
The interaction between the stiff particles and the matrix under large deformations highlights the significant complexity of the deformation patterns, especially under plane stress constraints.

\begin{figure}[ht]
	\begin{center}
		\begin{tabular}{c c}
            plane stress compressible & plane stress nearly incompressible \\
		    \includegraphics[scale=0.22]{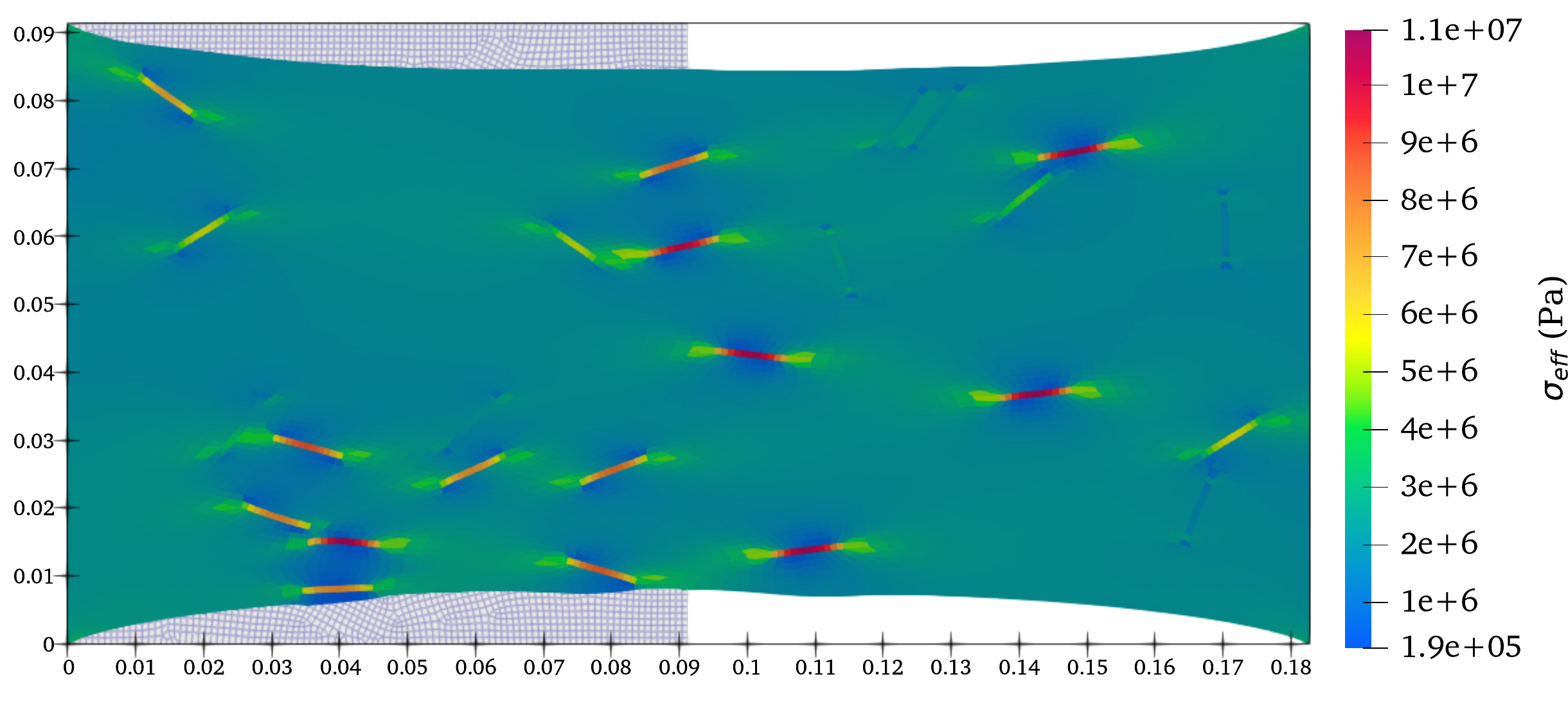} & \includegraphics[scale=0.22]{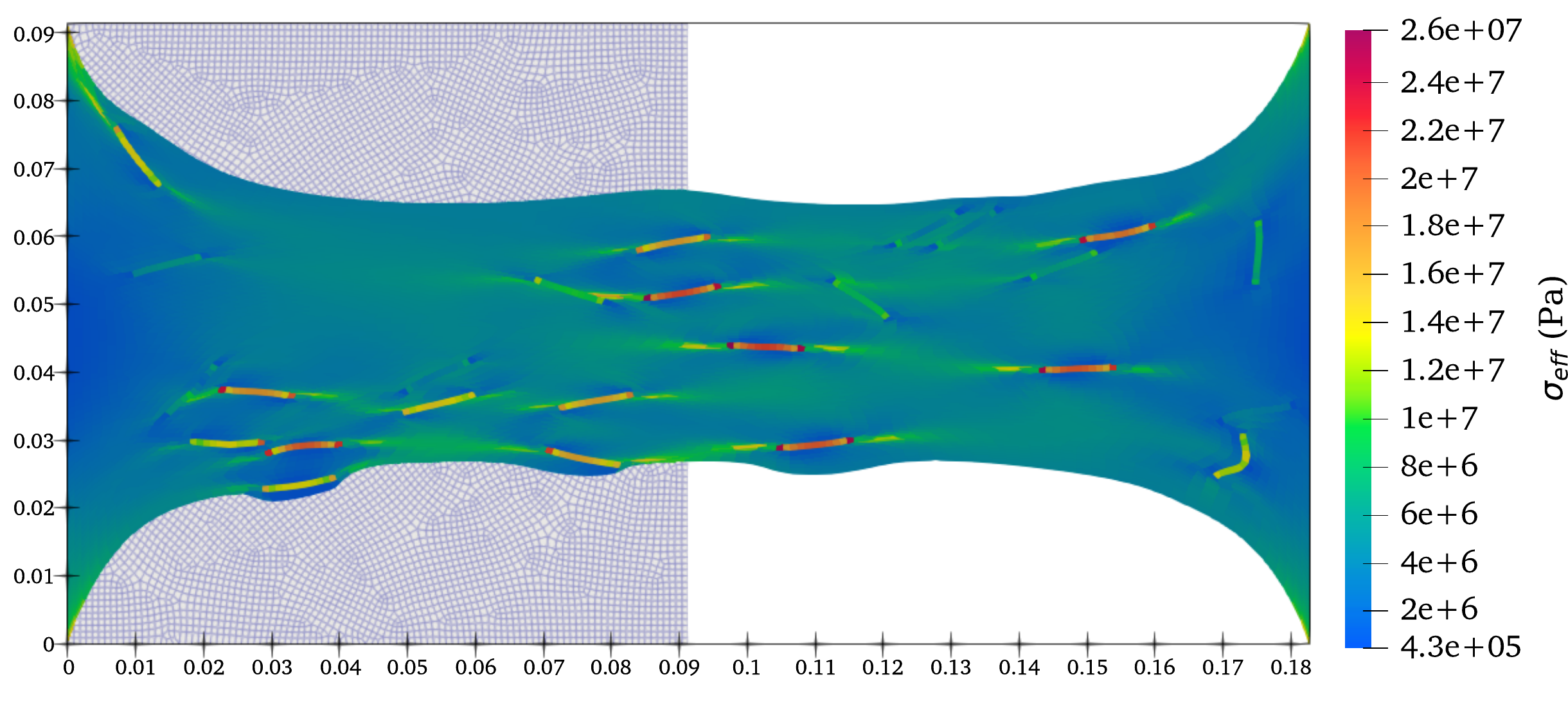}
		\end{tabular}
	\end{center}
    \caption{The plane stress compressible ($\nu=0.3$) and nearly incompressible ($\nu=0.4999$) stretched composite with short fibres.}
    \label{fig:pwf_contours}
\end{figure}

For the second scenario, a total of 25 rectangular fibres are distributed throughout a square block, with the volume fraction of the fibres set at $3 \%$.
The fibres are randomly dispersed in terms of both position and angle.
The aspect ratio of the fibres, defined as the ratio of fibre length to diameter, is chosen as $l/d = 10$.
Figure \ref{fig:pwf_contours} presents contour plots illustrating the von--Mises stress distribution within the square block, where the presence of embedded fibres significantly influences stress field.
The contour plots convey insights similar to those in the previous example, yet with a notable distinction.
Unlike the circular particles, the fibres undergo considerable bending and twisting, deviating from their initial shapes.
This behaviour arises from the elongated shape of fibres, contrasting with the uniform shape of circles despite both having equivalent stiffness.
The bending and twisting of fibres are more pronounced in the case of incompressible materials, primarily due to the matrix undergoing significantly more necking.
Moreover, it can be seen that the fibres aligned with the loading direction exhibit higher stress levels compared to those that are not aligned.
The extreme deformation of the fibres under stretching conditions emphasizes the complexity of the material response, showcasing the capabilities of the computational framework in capturing such detailed behaviour.

\section{Conclusion}\label{sec: conclusion}

A general formulation for the finite element modelling of the nonlinear deformation of arbitrary compressible (and nearly incompressible) hyperelastic materials under plane stress conditions along with its finite element code have been presented.
Additional insight into other common two-dimensional approaches has been provided.
The resulting numerical implementation is robust, demonstrating quadratic convergence for challenging problems.
The open-source finite element code is freely available online for broader use and adaptation \citep{madeal2024}.
Future research will extend this model to inelastic materials and focus on numerical aspects such as the choice of preconditioner and solver.

\section*{Acknowledgements}
Masoud Ahmadi's work is supported by the Engineering and Physical Sciences Research Council's Doctoral training fund.
Prashant Saxena acknowledges the financial support of the Engineering and Physical Sciences Research Council via project no. EP/V030833/1.
Paul Steinmann acknowledges the support of the European Research Council (Grant-No. 101052785, project: SoftFrac).

\section*{Supplementary data}\label{sec:sup}
The accompanying finite element codes are available open source at \url{https://doi.org/10.5281/zenodo.11636988}.

\bibliographystyle{myunsrtnat}
\bibliography{mybibfile}

\newpage
\appendix


\section{One-field finite element formulation}\label{sec: d and l}

\subsection{Variational setting}\label{sec:dirich}

Consider a hyperelastic body in the reference configuration subject to traction $\overline{\mathbf{T}}$ and body forces $\overline{\mathbf{B}}$,
with the whole boundary $\Gamma$ decomposed into Dirichlet $\Gamma^u$ and Neumann $\Gamma^t$ parts such that  $\Gamma = \Gamma^u \cup \Gamma^t$ and $\Gamma^u \cap \Gamma^t = \emptyset \, $.
The external potential energy generated is the negative of the work done by the applied loading, that is,
\begin{equation}
    \Pi^{\text{ext}} (\mathbf{u}) = -\int_{\Omega_X}{\mathbf{u} \cdot \overline{\mathbf{B}} \, \text{d} V}-\int_{\Gamma^t}{\mathbf{u}\cdot\overline{\mathbf{T}}}\, \text{d} A \, .
\end{equation}
The internal potential energy stored in the system due to deformation is given as 
\begin{equation}
    \Pi^{\text{int}} (\mathbf{u}) = \int_{\Omega_X}{\psi\left(\mathbf{F}\right) \text{d} V}.   
\end{equation}


\noindent
Thus, the total potential energy of the system expands as
\begin{equation}
    \Pi (\mathbf{u}) = \Pi^{\text{int}} (\mathbf{u}) + \Pi^{\text{ext}} (\mathbf{u}) .   
\end{equation}
Then, the variational formulation of the problem reads as
\begin{equation}
    \delta\Pi\left(\mathbf{u}; \delta\mathbf{u}\right) = \delta\Pi^{\text{int}}\left(\mathbf{u}; \delta\mathbf{u}\right) + \delta\Pi^{\text{ext}}\left(\mathbf{u}; \delta\mathbf{u}\right) = 0 \, , \quad \forall \delta \mathbf{u} \; \text{with} \; \delta \mathbf{u} = \mathbf{0} \; \text{on} \; \Gamma^u.
\end{equation}

\subsection{Linearization}
To linearize $\delta \Pi$ around a state of finite deformation, assuming dead loads, we introduce
%
\begin{equation}
    \mathcal{L} := \delta\Pi^{\text{int}}\left(\mathbf{u}; \delta\mathbf{u}\right) = \int_{\Omega_X}{\delta\mathbf{E}:\mathbf{S} \, \text{d}V} = \int_{\Omega_X}{\bm{\nabla}_x \, \delta \mathbf{u}: \bm{\tau} \, \text{d}V} \, .
\end{equation}
Denoting the linearized form of an expression by $\Delta (\bullet)$, 
one obtains
\begin{equation}\label{eq:DeltaL}   
    \Delta \mathcal{L} = \int_{\Omega_X}{\left[\delta\mathbf{E}:\Delta\mathbf{S}+\mathbf{S}:\Delta \delta\mathbf{E}\right] \text{d}V} .
\end{equation}
Applying the chain rule results in
\begin{equation}
    \Delta\mathbf{S}=\frac{\partial\mathbf{S}}{\partial\mathbf{E}}:\Delta\mathbf{E}=\frac{\partial^2\psi}{\partial\mathbf{E} \, \partial\mathbf{E}}:\Delta\mathbf{E} \, ,
\end{equation}
and the incremental constitutive tensor is identified by $\mathbb{C}=\partial^2\psi / [ \partial\mathbf{E} \, \partial\mathbf{E} ] \,$.
Using $\Delta\mathbf{E}=\text{sym}\left(\mathbf{F}^T \,[\bm{\nabla}_X \, \Delta\mathbf{u}] \right)$ and $\Delta \delta\mathbf{E} = \text{sym} ( \left[ \bm{\nabla}_X \, \delta\mathbf{u} \right]^T \, \bm{\nabla}_X \, \Delta\mathbf{u})$, together with the symmetry of $\mathbf{S}$ and the minor symmetry of $\mathbb{C}$, Equation \eqref{eq:DeltaL} is rewritten as
\begin{subequations}
\begin{align}
    \Delta \mathcal{L} &= \int_{\Omega_X}{\left[\mathbf{F}^T \left[ \bm{\nabla}_X \, \delta\mathbf{u} \right]:\mathbb{C}:\mathbf{F}^T \left[ \bm{\nabla}_X \, \Delta\mathbf{u} \right] +\left[ \bm{\nabla}_X \, \delta\mathbf{u} : \bm{\nabla}_{X} \, \Delta\mathbf{u} \right] \, \mathbf{S} \right] \text{d}V} \label{eq: DL1} \\
    &= \int_{\Omega_X}{\left[\bm{\nabla}_x \, \delta\mathbf{u}:J \, \mathbb{c}:\bm{\nabla}_x \, \Delta\mathbf{u}+ \left[ \bm{\nabla}_{x} \, \delta\mathbf{u} : \bm{\nabla}_{x} \, \Delta\mathbf{u} \right] \, \bm{\tau}\right] \text{d}V} , \label{eq: DL2}
\end{align}
\end{subequations}
where the following relations have been used 
\begin{align}
[\bm{\nabla}_x \, \delta\mathbf{u} ] \mathbf{F}=\bm{\nabla}_X \, \delta\mathbf{u} && \text{and} && \mathbb{c}_{ijkl}=J^{-1} \, F_{iI} \, F_{jJ}  \, F_{kK} \, F_{lL} \, \mathbb{C}_{IJKL} \, .    
\end{align}

Decomposing the incremental constitutive tensor into volumetric and isochoric contributions is also advantageous for numerical implementations.
Refer to Appendix \ref{app: ic tensor split} for the corresponding details. 

\subsection{Finite element approximation}\label{sec: fem}

The variational formulation derived in Section \ref{sec:dirich} is cast into the fully discrete form \citep{wriggers2008nonlinear, bonet1997nonlinear}. 
This involves discretizing the governing equations.
%
%
%

The primary variable in a classical single-field approach is the displacement $\mathbf{u}$. 
Consider a spatial discretisation of the domain into non-overlapping elements. 
The nodal variables are interpolated using displacement vector-valued shape functions $\mathbf{N}^{I,u}$.
The various fields are approximated as
\begin{align}
    \mathbf{u} = \sum_{I=1}^{n_{\text{dof}}}{\mathbf{N}^{I,u} \, u^I}, &&
    \delta \mathbf{u} = \sum_{I=1}^{n_{\text{dof}}}{\mathbf{N}^{I,u} \, {\delta u}^I}, && \Delta \mathbf{u}=\sum_{I=1}^{n_{\text{dof}}}{\mathbf{N}^{I,u} \, {\Delta u}^I} ,
\end{align}
where $n_{\text{dof}}$ is the number of global degrees of freedom.

Approximating the function $\delta\Pi$ using 
a first-order Taylor expansion gives
\begin{equation}
    \delta\Pi\left(\mathbf{u}^{t+1};\delta\mathbf{u}\right) \approx \delta\Pi \left( \mathbf{u}^t;\delta\mathbf{u} \right) + \Delta \delta \Pi \left( \mathbf{u}^t,;\delta\mathbf{u}, \Delta\mathbf{u} \right) \doteq 0 \, ,
\end{equation}
where $\mathbf{u}^{t+1}=\mathbf{u}^t+\Delta \mathbf{u}$ and the superscript $t$ denotes the load step number.
The equation is solved iteratively within a load step using a Newton scheme.
The second term is approximated by  Equation \eqref{eq: DL1}, that is $\Delta \, \delta\Pi\left(\mathbf{u};\delta\mathbf{u}, \Delta\mathbf{u}\right) = \Delta \mathcal{L}$. 

The residual vector $\mathbf{R}$, and the  tangent stiffness matrix $\mathbf{K}$, 
that form the linear system of equations 
$\mathbf{K} \Delta \mathbf{u} = \mathbf{R}$,
are obtained respectively as 
\begin{subequations}
\begin{align}
    R^I =& -\int_{\Omega_X}{\mathbf{S} : \bm{\nabla}_X \, \mathbf{N}^{I,u} \, \text{d}V} + \int_{\Omega_X}{\mathbf{N}^{I,u} \cdot {\overline{\mathbf{B}}} \, \text{d}V} + \int_{\Gamma^t}{\mathbf{N}^{I,u} \cdot \overline{\mathbf{T}}} \, \text{d}A \, , \\
    K^{IJ}=&\int_{\Omega_X}{ \mathbf{F}^T \left[ \mathbf{\nabla}_X \, \mathbf{N}^{I,u} \right] : \mathbb{C} : \mathbf{F}^T \left[ \mathbf{\nabla}_X \, \mathbf{N}^{J,u} \right] \, \text{d}V} + \int_{\Omega_X}{\left[ \mathbf{\nabla}_X \, \mathbf{N}^{I,u} : \mathbf{\nabla}_X \, \mathbf{N}^{J,u} \right] \mathbf{S} \, \text{d}V} .
\end{align}
\end{subequations}

\section{Three-field mixed finite element formulation}\label{app: 3field}

For nearly incompressible material,
due to locking, the single-field finite element formulation performs sub-optimally.
Various finite element formulations such as the B-bar and F-bar \citep{reese2002equivalent, de1996design, elguedj2008f}, enhanced assumed strain (EAS) \citep{sze2004stabilized, simo1990class, ye2020large}, selective reduced integration (SRI) \citep{hung2009addressing, simo1991quasi}, mixed methods \citep{simo1992geometrically, shojaei2018compatible, boffi2013mixed, schroder2011new, wriggers2008nonlinear} and stabilization techniques \citep{franca1988new, duddu2012finite, reese2000stabilization, maniatty2002higher, puso2006stabilized} have been developed to address locking.
Among these, mixed formulations which incorporate additional field variables, such as stress or strain fields, alongside the displacement field are particularly attractive \citep{simo1992geometrically, weiss1996finite, srinivasan2008generalized, zdunek2016five}.
Here, to avoid locking in the flatland, plane strain, and three-dimensional case, a three-field mixed finite element approximation is used, as detailed below.

\subsection{Variational setting}

The three-field formulation introduces the independent pressure $\widetilde{p}$ and strain variable $\widetilde{J}$ in addition to the displacement $\mathbf{u}$, expanding the set of unknowns to $\mathbf{Q}=\{\mathbf{u}, \widetilde{p}, \widetilde{J}\}$.
Hence, the constraint $\widetilde{J} \rightarrow J \left( \mathbf{u} \right)$ has to be fulfilled and the independent pressure follows as $\widetilde{p}
\rightarrow \partial\psi_{\text{vol}} (\widetilde{J}) / \partial\widetilde{J} \,$.
In the Hu--Washizu variational principle, the strain energy function is defined by
\begin{equation}
    \overline{\psi} (\mathbf{Q}) = \psi_{\text{iso}}(\widehat{\mathbf{C}}\left(\mathbf{u}\right))+\psi_{\text{vol}}(\widetilde{J})+\widetilde{p}\left[J(\mathbf{u})-\widetilde{J}\right] .
\end{equation}
While the external energy is the same as the single-field formulation, the internal potential energy for the three-field Hu--Washizu formulation takes the form
\begin{equation}
    \Pi^{\text{int}}\left(\mathbf{Q}\right) = \int_{\Omega_X}{\overline{\psi} \, \text{d}V}=\int_{\Omega_X}{\left[\psi_{\text{iso}}+\psi_{\text{vol}}+\widetilde{p}\left[J-\widetilde{J}\right]\right] \text{d}V}.
\end{equation}

\noindent
The first variation thereof reads
\begin{subequations}
\begin{align}
    D_{\delta\mathbf{u}}\Pi\left(\mathbf{Q}\right)&= \int_{\Omega_X}{\bm{\nabla}_x \, \delta \mathbf{u}: \bm{\tau} \, \text{d}V}-\int_{\Omega_X}{\delta\mathbf{u} \cdot \overline{\mathbf{B}} \, \text{d}V}-\int_{\Gamma^t}{\delta\mathbf{u} \cdot \overline{\mathbf{T}}} \, \text{d}A = 0 \, , \\
    D_{\delta \widetilde{p}}\Pi\left(\mathbf{Q}\right)&=\int_{\Omega_X}{\delta \widetilde{p} \left[J(\mathbf{u})-\widetilde{J}\right] \text{d}V}=0 \, , \\
    D_{\delta\widetilde{J}}\Pi\left(\mathbf{Q}\right)&=\int_{\Omega_X}{\delta\widetilde{J} \left[\frac{\text{d}\psi_{\text{vol}}(\widetilde{J})}{\text{d}\widetilde{J}}-\widetilde{p}\right]}\text{d}V=0 \, .
\end{align}
\end{subequations}

\subsection{Linearization}
Applying the approach that led to \eqref{eq: DL2}, renders for the three-field problem:
\begin{equation}
    \Delta\delta\Pi\left(\mathbf{Q}; \delta\mathbf{Q}, \Delta\mathbf{Q}\right) \Rightarrow \begin{cases}
        D_{\Delta \mathbf{u}} \delta \Pi (\mathbf{Q}; \delta \mathbf{Q}) =\int_{\Omega_X}{ [\bm{\nabla}_x \, \delta \mathbf{u}:J\mathbb{c}+\bm{\nabla}_x \, \delta \mathbf{u} \, \bm{\tau} +\delta \widetilde{p} J \mathbf{I} : \bm{\nabla}_x \, \Delta \mathbf{u} ] \text{d}V} \\
        D_{\Delta \widetilde{p}} \delta \Pi (\mathbf{Q}; \delta \mathbf{Q}) = \int_{\Omega_X}{ \left[ \bm{\nabla}_x \, \delta \mathbf{u}:J \mathbf{I}-\delta \widetilde{J} \Delta \widetilde{p} \, \right] \text{d}V} \\
        D_{\Delta \widetilde{J}} \delta \Pi (\mathbf{Q}; \delta \mathbf{Q}) = \int_{\Omega_X}{\left[ -\delta \widetilde{p} + \delta \widetilde{J} \, \dfrac{\text{d}^2\psi_{\text{vol}}}{\text{d}\widetilde{J} \, \text{d}\widetilde{J}} \right]\Delta \widetilde{J} \, \text{d}V}
    \end{cases} \, .
\end{equation}
Approximating the function $\delta\Pi$ using a first-order Taylor expansion gives
\begin{equation}
    \delta\Pi\left(\mathbf{Q}^{t+1};\delta\mathbf{Q}\right) \approx \delta\Pi\left(\mathbf{Q}^t;\delta\mathbf{Q}\right) + \Delta \, \delta\Pi\left(\mathbf{Q}^t; \delta\mathbf{Q}, \Delta\mathbf{Q}\right) = 0 \, ,
\end{equation}
and by setting $\mathbf{Q}^{t+1}=\mathbf{Q}^t+\Delta \mathbf{Q}$, one can solve the system of equations iteratively employing a Newton scheme.

\subsection{Finite element approximation}
Denoting the nodal variables at global node $I$ as $\mathbf{Q}^I=\left\{\mathbf{u}^I,\widetilde{p}^I,{\widetilde{J}}^I\right\}$, one can approximate the various fields using basis functions as
\begin{align}
    \mathbf{u} = \sum_{I \in \mathcal{I}_u}{\mathbf{N}^{I,u} \, u^I}, &&\widetilde{p}=\sum_{I \in \mathcal{I}_{\widetilde{p}}}{N^{I,\widetilde{p}} \, \widetilde{p}^I}, && \widetilde{J}=\sum_{I \in \mathcal{I}_{\widetilde{J}}}{N^{I,\widetilde{J}} \, {\widetilde{J}}^I} \, .
\end{align}
The sets $\mathcal{I}_u$, $\mathcal{I}_{\widetilde{p}}$ and $\mathcal{I}_{\widetilde{J}}$ contain the global degrees of freedom for the respective variables.
Similarly, the components of $\delta\mathbf{Q}$ and $\Delta\mathbf{Q}$ can also be approximated as
\begin{align}
    \delta \mathbf{u}&=\sum_{I \in \mathcal{I}_u}{\mathbf{N}^{I,u} \, {\delta u}^I}, \quad \delta \widetilde{p}=\sum_{I \in \mathcal{I}_{\widetilde{p}}}{N^{I,\widetilde{p}} \, {\delta \widetilde{p}}^I}, \quad \delta\widehat{J}=\sum_{I \in \mathcal{I}_{\widetilde{J}}}{N^{I,\widetilde{J}} \, {\delta\widetilde{J}}^I}, \nonumber \\
    \Delta \mathbf{u}&=\sum_{I \in \mathcal{I}_u}{\mathbf{N}^{I,u} \, {\Delta u}^I}, \quad \Delta \widetilde{p}=\sum_{I \in \mathcal{I}_{\widetilde{p}}}{N^{I,\widetilde{p}} \, {\Delta \widetilde{p}}^I}, \quad \Delta\widetilde{J}=\sum_{I \in \mathcal{I}_{\widetilde{J}}}{N^{I,\widehat{J}} \, \Delta{\widetilde{J}}^I}.
\end{align}
Then the residual vector $\mathbf{R} = \left\{\mathbf{R}_u, \mathbf{R}_{\widetilde{p}}, \mathbf{R}_{\widetilde{J}}\right\}$ follows as
\begin{subequations}
\begin{align}    
    R_u^I& = -\int_{\Omega_X}{\bm{\tau}: \bm{\nabla}_x \, \mathbf{N}^{I,u} \, \text{d}V} + \int_{\Omega_X}{\mathbf{N}^{I,u} \cdot \overline{\mathbf{B}} \, \text{d}V} + \int_{\Gamma^t}{\mathbf{N}^{I,u} \cdot  \overline{\mathbf{T}}} \, \text{d}A, \\
    R_{\widetilde{p}}^I& = -\int_{\Omega_X}{N^{I,\widetilde{p}} \left[J-\widetilde{J}\right] \, \text{d}V}, \\
    R_{\widetilde{J}}^I& = -\int_{\Omega_X}{N^{I,\widetilde{J}}\left[\frac{\text{d}\psi_{\text{vol}}(\widetilde{J})}{\text{d}\widetilde{J}}-\widetilde{p}\right]}\text{d}V .
\end{align}
\end{subequations}
Finally, the non-zero components of the tangent stiffness matrix $\mathbf{K}$ associated with the degrees of freedom $I, J \in \{\mathcal{I}_u, \mathcal{I}_{\widetilde{p}}, \mathcal{I}_{\widetilde{J}}\}$ read
\begin{subequations}
\begin{align}
    K_{uu}^{IJ}& = \int_{\Omega_X}{\left[ \bm{\nabla}_x \, \mathbf{N}^{I,u} : J \mathbb{c} : \bm{\nabla}_x \, \mathbf{N}^{J,u} + \left[ \bm{\nabla}_x \, \mathbf{N}^{I,u} : \bm{\nabla}_x \, \mathbf{N}^{J,u} \right] \bm{\tau} \right] \, \text{d}V}, \quad &\forall I \in \mathcal{I}_u  \, \text{and} \, \forall J \in \mathcal{I}_u \, , \\
    K_{u \widetilde{p}}^{IJ}&=\int_{\Omega_X}{\bm{\nabla}_x \, \mathbf{N}^{I,u} : J \, \mathbf{I}} \, N^{J,\widetilde{p}} \, \text{d}V, \quad &\forall I \in \mathcal{I}_u  \, \text{and} \, \forall J \in \mathcal{I}_{\widetilde{p}} \, , \\
    K_{\widetilde{p} u}^{IJ}&=\int_{\Omega_X}{N^{I,\widetilde{p}}\, J\, \mathbf{I} : \bm{\nabla}_x \, N_k^{J,u} \, \text{d}V}, \quad &\forall I \in \mathcal{I}_{\widetilde{p}}  \, \text{and} \, \forall J \in \mathcal{I}_u \, , \\
    K_{\widetilde{p}\widetilde{J}}^{IJ}&=\int_{\Omega_X}{-N^{I,\widetilde{p}}\, N^{J,\widehat{J}}} \, \text{d}V, \quad &\forall I \in \mathcal{I}_{\widetilde{p}}  \, \text{and} \, \forall J \in \mathcal{I}_{\widetilde{J}} \, , \\
    K_{\widetilde{J}\widetilde{p}}^{IJ}&=\int_{\Omega_X}{-N^{I,\widehat{J}}}N^{J,\widetilde{p}} \, \text{d}V, \quad &\forall I \in \mathcal{I}_{\widetilde{J}}  \, \text{and} \, \forall J \in \mathcal{I}_{\widetilde{p}} \, , \\
    K_{\widetilde{J}\widetilde{J}}^{IJ}&=\int_{\Omega_X}{N^{I,\widehat{J}} \left[\frac{d^2\psi_{\text{vol}}}{\text{d}\widetilde{J} \, \text{d}\widetilde{J}}\right]N^{J,\widehat{J}}}\text{d}V, \quad &\forall I \in \mathcal{I}_{\widetilde{J}}  \, \text{and} \, \forall J \in \mathcal{I}_{\widetilde{J}} \, ,
\end{align}
\end{subequations}
where
\begin{equation}
    \mathbf{K}=\left[\begin{matrix}\mathbf{K}_{uu}&\mathbf{K}_{u \widetilde{p}}&\mathbf{K}_{u \widetilde{J}} \\ 
    \mathbf{K}_{\widetilde{p}u}&\mathbf{K}_{\widetilde{p}\widetilde{p}}&\mathbf{K}_{\widetilde{p}\widetilde{J}} \\ 
    \mathbf{K}_{\widetilde{J} u}&\mathbf{K}_{\widetilde{J}\widetilde{p}}&\mathbf{K}_{\widetilde{J}\widetilde{J}}\\\end{matrix}\right] =\left[\begin{matrix}\mathbf{K}_{uu}&\mathbf{K}_{u \widetilde{p}}& \bm{0} \\
    \mathbf{K}_{\widetilde{p} u}&\bm{0}&\mathbf{K}_{\widetilde{p}\widetilde{J}}\\
    \bm{0}&\mathbf{K}_{\widetilde{J}\widetilde{p}}&\mathbf{K}_{\widetilde{J}\widetilde{J}}\\\end{matrix}\right] .
\end{equation}

\section{Isochoric-volumetric decomposition}

\subsection{Decomposition of the stress} \label{app: stress split}

Based on \eqref{eq: psi split}, the Piola--Kirchhoff stress tensor, $\mathbf{S}$, can be stated as
\begin{equation}\label{eq:Ssplit}
    \mathbf{S} = 2 \, \frac{\partial\psi(\widehat{\mathbf{C}},J)}{\partial\mathbf{C}}=2\left[\frac{\partial\psi_{\text{iso}}(\widehat{\mathbf{C}})}{\partial\mathbf{C}}+\frac{\partial\psi_{\text{vol}}\left(J\right)}{\partial\mathbf{C}}\right] .
\end{equation}
By using the chain rule of differentiation on the first term, one obtains
\begin{equation}
    \frac{\partial\psi_{\text{iso}}}{\partial\mathbf{C}}=\frac{\partial\psi_{\text{iso}}}{\partial\widehat{\mathbf{C}}} : \frac{\partial\widehat{\mathbf{C}}}{\partial\mathbf{C}} \, ,
\end{equation}
where
\begin{equation}
    \frac{\partial\widehat{\mathbf{C}}}{\partial\mathbf{C}} = \frac{\partial\left(J^{-2/3} \, \mathbf{C}\right)}{\partial\mathbf{C}} =\frac{\partial J^{-2/3}}{\partial\mathbf{C}}\otimes\mathbf{C}+J^{-2/3} \, \frac{\partial\mathbf{C}}{\partial\mathbf{C}}=J^{-2/3}\left[\mathbb{I}-\frac{1}{3} \, \mathbf{C}^{-1}\otimes\mathbf{C}\right] \equiv \mathbb{P} \, .
\end{equation}
For the second part of Eq. \eqref{eq:Ssplit}, the chain rule is used again to obtain
\begin{equation}
    \frac{\partial\psi_{\text{vol}}}{\partial\mathbf{C}}=\frac{\partial\psi_{\text{vol}}}{\partial J}\frac{\partial J}{\partial\mathbf{C}} \, ,
\end{equation}
where $\dfrac{\partial J}{\partial \mathbf{C}}=\frac{1}{2} J \mathbf{C}^{-1}$ and $p = \dfrac{ \partial \psi_{\text{vol}} (J)}{\partial J}$ is recognised as the hydrostatic pressure.
Hence,
\begin{equation}
    \mathbf{S}= \mathbb{P} \left[2 \, \frac{\partial\psi_{\text{iso}}}{\partial\widehat{\mathbf{C}}}\right]+pJ \, \mathbf{C}^{-1}=\mathbf{S}_{\text{iso}}+\mathbf{S}_{\text{vol}} \, .
\end{equation}
Transforming this relation to the current configuration, the Kirchhoff stress can be written as
\begin{equation}
    \bm{\tau}=\mathbf{F} \, \mathbf{S} \, \mathbf{F}^T=pJ \, \mathbf{I} + \mathbb{D}:\widehat{\bm{\tau}}=\bm{\tau}_{\text{vol}}+\bm{\tau}_{\text{iso}} \, ,
\end{equation}
where 
$\mathbb{D}=\mathbb{I}-\frac{1}{3} \, \mathbf{I}\otimes\mathbf{I}$, and $\widehat{\bm{\tau}}$ (for isotropy) follows as
\begin{equation}
    \widehat{\bm{\tau}}=\widehat{\mathbf{F}} \left[2 \, \frac{\partial\psi_{\text{iso}}(\widehat{\mathbf{C}})}{\partial\widehat{\mathbf{C}}}\right]{\widehat{\mathbf{F}}}^T=2\widehat{\mathbf{b}} \, \frac{\partial\psi_{\text{iso}}(\widehat{\mathbf{b}})}{\partial\widehat{\mathbf{b}}} \, .
\end{equation}

\subsection{Decomposition of the incremental constitutive tensor}\label{app: ic tensor split}

The fourth-order incremental constitutive tensor in Lagrangian description is decomposed into isochoric and volumetric parts as follows
\begin{equation}
    \mathbb{C}=\mathbb{C}_{\text{iso}}+\mathbb{C}_{\text{vol}} \, .
\end{equation}
For isotropy, the incremental constitutive tensor in the Eulerian description is defined by 
\begin{equation}
    \mathbb{c}=4 J^{-1} \, \mathbf{b} \, \frac{\partial^2\psi}{\partial\mathbf{b} \, \partial\mathbf{b}} \, \mathbf{b} \, 
\end{equation}
and relates to its Lagrangian counterpart by the push-forward
\begin{equation}
    c_{ijkl}=J^{-1} \, F_{iI}F_{jJ} \, C_{IJKL} \, F_{kK}F_{lL} \, .
\end{equation}
It likewise decomposes into isochoric and volumetric parts as
\begin{equation}
    \mathbb{c}=\mathbb{c}_{\text{iso}}+\mathbb{c}_{\text{vol}} \,  ,
\end{equation}
where
\begin{equation}
    \mathbb{c}_{\text{vol}}=4 J^{-1} \, \mathbf{b} \, \frac{\partial^2\psi_{\text{vol}}}{\partial\mathbf{b} \, \partial\mathbf{b}} \, \mathbf{b}=\left[p+J \, \frac{dp}{dJ}\right]\mathbf{I} \otimes \mathbf{I}-2p\, \mathbb{I} \, ,
\end{equation}
and
\begin{equation}
    J \, \mathbb{c}_{\text{iso}}=4\mathbf{b} \, \frac{\partial^2\psi_{\text{iso}}}{\partial\mathbf{b} \, \partial\mathbf{b}} \, \mathbf{b}=\mathbb{D}:\widehat{\mathbb{c}}:\mathbb{D}+\frac{2}{3} \text{tr} (\widehat{\bm{\tau}}) \, \mathbb{D}-\frac{2}{3}\left[\mathbf{I}\otimes\bm{\tau}_{\text{iso}}+\bm{\tau}_{\text{iso}}\otimes\mathbf{I}\right] .
\end{equation}
Here the abbreviation $\widehat{\mathbb{c}}$ is  defined by
\begin{equation}
    \widehat{\mathbb{c}}=4\widehat{\mathbf{b}} \, \frac{\partial^2\psi_{\text{iso}}}{\partial\widehat{\mathbf{b}} \, \partial\widehat{\mathbf{b}}} \, \widehat{\mathbf{b}} \, . 
\end{equation}

\end{document}